\documentclass[oneside,final]{amsart}
\usepackage{geometry}

\usepackage[utf8]{inputenc}
\usepackage{hyperref}
\usepackage[english]{babel}
\usepackage{enumerate}
\usepackage{geometry}
\usepackage{aliascnt}
\usepackage{amssymb}
\usepackage{amsmath}
\usepackage{amsthm}
\usepackage{verbatim}
\usepackage[normalem]{ulem}
\usepackage{mathtools}
\usepackage{esint} 
\usepackage{nth}
\usepackage{ragged2e}

\usepackage[noabbrev,capitalize]{cleveref}

\usepackage[right]{showlabels}

\usepackage{tikz}
\usepackage{tikzscale}
\usepackage{graphicx}
\usepackage[labelformat=simple]{subcaption}
\usepackage{mathscinet} 

\usepackage[title]{appendix}

\usepackage{color}

\numberwithin{equation}{section} 

\def\todo{\marginpar{\textcolor{red}{todo}}\textcolor{red}}

\usepackage{dsfont} 
\usepackage{hyperref} 

\usepackage{autonum} 

\allowdisplaybreaks[4] 

\theoremstyle{plain}

\newtheorem{proposition}{Proposition}[section]

\newaliascnt{lemma}{proposition} 

\newtheorem{lemma}[lemma]{Lemma}

\aliascntresetthe{lemma} 
\Crefname{lemma}{Lemma}{Lemmas}

\newaliascnt{theorem}{proposition} 
\newtheorem{theorem}[theorem]{Theorem}

\aliascntresetthe{theorem} 

\newaliascnt{corollary}{proposition} 
\newtheorem{corollary}[corollary]{Corollary}

\aliascntresetthe{corollary}

\newaliascnt{hypothesis}{proposition}

\aliascntresetthe{hypothesis} 

\theoremstyle{definition}

\newaliascnt{definition}{proposition} 
\newtheorem{definition}[definition]{Definition}

\aliascntresetthe{definition} 
\Crefname{definition}{Definition}{Definitions}

\newaliascnt{problem}{proposition} 

\aliascntresetthe{problem} 

\newaliascnt{example}{proposition} 

\aliascntresetthe{example} 

\theoremstyle{remark}

\newaliascnt{remark}{proposition} 
\newtheorem{remark}[remark]{Remark}

\aliascntresetthe{remark} 

\def\equationautorefname~#1\null{%
	(#1)\null
}
 
\newcommand{\R}{\mathbb{R}}

\newcommand{\N}{\mathbb{N}}

\newcommand{\V}{\mathcal{V}}

\newcommand{\scurv}{\kappa}
\renewcommand{\S}{\mathbb{S}}

\newcommand{\A}{\mathcal{A}}

\renewcommand{\H}{\mathcal{H}}
\newcommand{\W}{\mathcal{W}}

\newcommand{\Ll}{\mathcal{L}}

\newcommand{\genus}{\operatorname{genus}}
\newcommand{\diam}{\mathrm{diam}}

\newcommand*{\dd}{\mathop{}\!\mathrm{d}}

\def\nicefrac#1#2{%
    \raise.5ex\hbox{$#1$}%
    \kern-.15em/\kern-.05em%
    \lower.25ex\hbox{$#2$}}

\title{Gradient flow dynamics for cell membranes in the Canham--Helfrich model}

\author[F.~Rupp]{Fabian Rupp}
\address[F.~Rupp]{Faculty of Mathematics, University of Vienna, Oskar-Morgenstern-Platz 1, 1090 Vienna, Austria.}
\email{fabian.rupp@univie.ac.at}
\author[C.~Scharrer]{Christian Scharrer}
\address[C.~Scharrer]{Institute for Applied Mathematics, University of Bonn, Endenicher Allee 60, 53115 Bonn, Germany.}
\email{scharrer@iam.uni-bonn.de}
\author[M.~Schlierf]{Manuel Schlierf}
\address[M.~Schlierf]{Salzburg University, Hellbrunner Str.~34, 5020 Salzburg, Austria.}
\email{math@manuelschlierf.info}

\begin{document}

\begin{abstract}
    The energetically most efficient way how a deformed red blood cell regains equilibrium is mathematically described by the gradient flow of the Canham--Helfrich functional, including a spontaneous curvature and the conservation of surface area and enclosed volume. Using a recently 
    discovered multiplicity inequality, we prove global existence and convergence of smooth solutions for spheres and axisymmetric tori, provided the initial energy lies below explicit thresholds.
\end{abstract}
\maketitle



\section{Introduction}

We consider the bending energy of an immersed surface $f\colon \Sigma \to \R^3$ defined by
\begin{equation}\label{eq:W}
	\mathcal W(f)\vcentcolon = \frac14\int_{\Sigma}H^2\,\mathrm d\mu
\end{equation}
where $H \vcentcolon= \kappa_1+\kappa_2$ is the mean curvature, and $\mu$ is the Riemannian measure on $\Sigma$ induced by the pull back metric $g_f$. This energy functional is not only \emph{scale-invariant,} but even invariant under all conformal transformations of the ambient space, as already observed by Thomsen \cite{Thomsen1923} and Blaschke \cite{Blaschke1929}. The energy $\mathcal{W}$ became more popular after the work of Willmore \cite{Willmore65} who also showed that among all immersions of closed surfaces, the round spheres are the absolute minimizers with energy $4\pi$ \cite{Willmore93Oxford}. An essential tool in examining the Willmore energy is an inequality of Li and Yau \cite{liyau1982} relating the energy with the multiplicity at any point. That is, if $p\in \R^3$ with $f^{-1}(p)=\{x_1,\ldots,x_k\}$, $x_i\neq x_j$ for $i\neq j$, then
\begin{equation} \label{eq:intro:LY}
    4\pi k \le \mathcal W(f).
\end{equation}
In particular, the condition $\mathcal{W}(f)<8\pi$ rules out points of higher multiplicity, implying that $f$ is an embedding. Using the conformal invariance, Bryant \cite{bryant1984} proved that \emph{Willmore spheres}, i.e., spheres critical for \eqref{eq:W}, always attain equality in \eqref{eq:intro:LY} and arise from the stereographic projection of minimal surfaces in $\mathbb{S}^3$. In particular, he proved that the only Willmore spheres with energy strictly less than $8\pi$ are the round spheres. This fact should become a key ingredient in the works of Kuwert and Schätzle on the \emph{Willmore flow}, the $L^2$-gradient flow of the Willmore functional \cite{kuwertschaetzle2001,kuwertschaetzle2002,kuwertschaetzle2004}. Unless additional symmetry is assumed \cite{dallacquamullerschatzlespener2020,Jakob2023}, the only generic long time existence result is for topological spheres due to the lack of an analogue of Bryant's classification result.

\subsection{The stationary Canham--Helfrich model}
One of the most prominent applications of the Willmore energy is the connection to biological membranes. In the 70s, Helfrich \cite{helfrich1973} derived the curvature elastic energy per unit area of a lipid bilayer
\begin{equation}\label{eq:intro:Helfrich}
    \frac{1}{2}k_c(H-c_0)^2 + \bar k_cK
\end{equation}
where $k_c,\bar k_c$ are curvature elastic moduli, see also the previous work of Canham \cite{canham1970}. The so-called spontaneous curvature $c_0\in\R$ is a parameter that allows the model to be adjusted with respect to both, the type of membrane as well as its two aqueous solutions. For the vital example of red blood cells, this parameter was found to be negative and large in modulus \cite{DeulingHelfrich}. The membrane itself is modeled as an immersed surface $f\colon \Sigma\to\R^3$ whose total stored energy is obtained by integrating \eqref{eq:intro:Helfrich} over $\Sigma$, where the Gauss--Bonnet theorem implies that the integral of $\bar k_cK$ is a topological constant.
Introducing $c_0\neq 0$ reduces the invariance group of the functional compared to the Willmore functional \eqref{eq:W}, see \cite{DeMatteisManno14}, so that previous methods based on conformal invariance \cite{liyau1982,bryant1984} are no longer available. Moreover, the energy is invariant only under \emph{positive} (i.e., orientation-preserving) reparametrizations.
Additional complexity arises from the inextensibility of the membrane and the osmotic pressure balance between the cell and its surroundings that are mathematically modeled by prescribing the area and volume. Given $c_0\in\R$, $g\in\N_0$, $A_0>0$, and $0<V_0 < A_0^{3/2}/\sqrt{36\pi}$, we thus arrive at the following minimization problem. 
\begin{align}
\label{intro:HP}\tag{HP}
\begin{split}
    &\mathrm{Minimize}\qquad \mathcal H_{c_0}(f)\vcentcolon=\frac{1}{4}\int_\Sigma(H-c_0)^2\,\mathrm d\mu, \\ & \text{among smooth immersions $f\colon\Sigma\to \R^3$ with } \text{$\operatorname{genus}\Sigma=g$, $\mathcal A(f)=A_0$, and $\mathcal V(f)=V_0$},
\end{split}
\end{align}
where $\mathcal{A}$ and $\mathcal{V}$ denote the area and volume.

Existence of solutions for $c_0=0$, $g=0$ was proven by Schygulla \cite{Schygulla}. For nonzero spontaneous curvature, Mondino and the second author \cite{MondinoScharrer20} showed that solutions do exist in general, but are not always given by smooth submanifolds. Indeed, a bubbling phenomenon may occur, cf. \cite[Example 1.2]{ruppscharrer2023}. For $c_0=0$ and $g\in\N_0$ arbitrary, existence of smooth minimizers has been proven recently, see \cite{KellerMondinoRiviere2014,Scharrer22NLA,MondinoScharrer2023,KusnerMcGrath23}, whereas existence in weak settings (without prescribing the topology) has been obtained in \cite{BrazdaLussardiStefanelli20,KubinLussardiMorandotti24}.
An explicit general criterion on the data $A_0,V_0,c_0$ that guarantees smoothly embedded spherical minimizers has been obtained by the first two authors \cite{ruppscharrer2023} by generalizing the inequality of Li and Yau \eqref{eq:intro:LY} for the Helfrich functional
\begin{equation}\label{eq:intro:HLY}
    4k\pi \le \mathcal H_{c_0}(f) - 2c_0\int_\Sigma\frac{\langle f - p,\nu\rangle}{|f-p|^2}\,\mathrm d\mu
\end{equation}
where $\nu$ is the unit normal induced by the orientation of $\Sigma$, cf.\ \Cref{sec:def}.
If $f$ is an Aleksandrov immersion (see \cite[Definition 1.4]{ruppscharrer2023}) or a volume varifold (see \cite[Hypothesis 4.5]{ruppscharrer2023}, \cite[Definition 3.1]{Scharrer23}), the last integral can be controlled in terms of the energy $\mathcal H_{c_0}$ and the data $A_0,V_0,c_0$, by means of the diameter bounds in \cite{ruppscharrer2023,Scharrer23}. Consequently, smooth solutions for the Helfrich problem \eqref{intro:HP} do exist, provided the minimal energy satisfies
\begin{equation}\label{intro:eq:energy-condition}
    \inf_{f}\mathcal H_{c_0}(f) < \Omega(A_0,V_0,c_0),
\end{equation}
where $\Omega(A_0,V_0,c_0)$ is defined as in \Cref{def:Omega} below.

\subsection{Dynamic approach and main results}
Capillaries, the smallest human blood vessels, can be so thin that red blood cells have to be deformed in order to fit through. Classical experiments show how red blood cells deform when sucked into a micropipette \cite{RandBurton}. 
After being released, the cells regain their original shape within seconds. 
In a  viscous regime, the continuous deformation
towards an energetically more favorable state can be schematically described by the $L^2$-gradient flow of the energy. 
Given an initial 
immersion $f_0\colon\mathbb S^2\to\R^3$,
a family of immersions $f\colon[0,T)\times \mathbb S^2\to \R^3$ with 
$f(0,\cdot)=f_0$ evolves according to the $L^2$-gradient flow of the Helfrich functional if
\begin{align}\label{eq:intro:proj_gradient-flow}
    \partial_t f = - 2\Pi(f) \nabla\mathcal{H}_{c_0}(f),
\end{align}
where $\Pi(f)$ is the $L^2(\dd \mu_f)$-orthogonal projection corresponding to the tangent space of the manifold of immersions with prescribed area and volume. This projection makes \eqref{eq:intro:proj_gradient-flow} a quasilinear nonlocal parabolic system, which is degenerate due to its geometric invariance. Moreover, being of fourth order, maximum principles are not available.
Whenever the mean curvature is not constant on $\Sigma$, i.e., the \emph{CMC-deficit} 
\begin{align}\label{eq:CMC-deficit}
    \overline{\mathcal{H}}(f) \vcentcolon = \frac{1}{4}\int_\Sigma(H-\overline{H})^2\dd \mu,
\end{align}
is strictly positive, 
Equation \eqref{eq:intro:proj_gradient-flow} can be equivalently written as
\begin{equation}\label{eq:intro:gradient-flow}
    \partial_t f = - \big( 2\nabla \mathcal H_{c_0}(f) + \lambda_1(f)\nabla \mathcal A(f) + \lambda_2(f)\nabla \mathcal V(f)\big)
\end{equation}
where the two Lagrange-multipliers $\lambda_1,\lambda_2$ are uniquely determined by the condition that area and volume are preserved along the flow, see \eqref{eq:def-l1} and \eqref{eq:def-l2} in \Cref{sec:preserving_Helfrich_flow}. Any smooth solution $f\colon [0,T)\times \Sigma \to\R^3$ of \eqref{eq:intro:gradient-flow} with constant area $A_0\equiv \mathcal A(f)$ and constant volume $V_0\equiv\mathcal V(f)$ will be referred to as an \emph{$(A_0,V_0,c_0)$-Helfrich flow}.
Short time existence and uniqueness has been proven by Nagasawa and Yi \cite{nagasawayi2012}.
In this paper, we prove global existence and convergence if the initial energy is below an explicit threshold determined by the parameters.

\begin{theorem}\label{thm:main-result}
    Consider $A_0,V_0>0$ with 
    \begin{equation}\label{eq:iso}
      \sigma\vcentcolon= 36\pi \frac{V_0^2}{A_0^3} \in (0,1),
    \end{equation}
    $c_0\in\R$, and let $f_0\colon\S^2\to\R^3$ be an embedding 
    with $\A(f_0)=A_0$ and $\V(f_0)=V_0$, satisfying
    \begin{equation}\label{eq:en-thresh}
        \sqrt{\H_{c_0}(f_0)} \leq \min\Bigl\{\sqrt{\Omega(A_0,V_0,c_0)}, \sqrt{\frac{4\pi}{\sigma}} - \frac{c_0}2\sqrt{A_0}\Bigr\}.
    \end{equation}

    Then, the $(A_0,V_0,c_0)$-Helfrich flow with initial datum $f_0$ exists globally, is embedded for all $t\geq0$, and converges after positive reparametrization smoothly to an embedded \emph{constrained $\mathcal{H}_{c_0}$-Helfrich sphere} $f_\infty\colon \mathbb{S}^2\to\R^3$ as $t\to\infty$, i.e., to a solution of
    \begin{align}\label{eq:intro:Helfrich_eq}
    \Delta H + |A^0|^2H +  c_0(\frac12H^2-|A^0|^2) - (\lambda_1+\frac12c_0^2) H - \lambda_2 = 0 \quad \text{ for some }\lambda_1,\lambda_2\in\R.
\end{align}
\end{theorem}

The assumption of spherical topology in \Cref{thm:main-result} is necessary in view of the absence of a classification result as in \cite{bryant1984} for higher genus Willmore surfaces.
However, if the initial surface is 
an axisymmetric torus, we may use the strategy from \cite{dallacquamullerschatzlespener2020} to conclude the following.
\begin{theorem}\label{intro:cor:axisymmetric_tori}
Let $A_0,V_0>0$ satisfy \eqref{eq:iso}, $c_0\in\R$, and $f_0\colon \mathbb{T}^2\to\R^3$ be an axisymmetric embedded torus with $\mathcal{A}(f_0)=A_0$ and $\mathcal{V}(f_0)=V_0$, satisfying
\eqref{eq:en-thresh}.
Then, the $(A_0,V_0,c_0)$-Helfrich flow with initial datum $f_0$ exists globally, is embedded for all $t\geq0$, and converges after positive reparametrization smoothly to an embedded axisymmetric constrained $\mathcal{H}_{c_0}$-Helfrich torus $f_\infty\colon \mathbb{T}^2\to\R^3$ as $t\to\infty$.
\end{theorem}

Whenever the parameters $A_0,V_0,c_0$ allow for an admissible initial datum for \Cref{thm:main-result,intro:cor:axisymmetric_tori}, we can conclude existence of critical points.
Since
\begin{equation}
    \lim_{c_0\to0}\mathcal H_{c_0} = \mathcal W,\qquad \lim_{c_0\to0} \Omega(A_0,V_0,c_0) = 8\pi,
\end{equation}
(see \Cref{lem:will-vs-helf} for a quantitative upper bound of the first limit) we have by \cite[Theorem 7.1]{rupp2024} that for all $A_0,V_0>0$ satisfying \eqref{eq:iso}, there exists $\varepsilon=\varepsilon(A_0,V_0)>0$ such that Condition \eqref{eq:en-thresh} admits axisymmetric spherical initial 
embeddings $f_0$ whenever $|c_0|<\varepsilon$. In addition, these $f_0$ are symmetric with respect to a plane orthogonal to the axis of rotation. 

On the other hand, since for $c_0<0$,
\begin{equation}
    \lim_{A_0\to0}\mathcal H_{c_0} = \mathcal W,\qquad \lim_{r\to 0+}\Omega(r^2A_0,r^3V_0,c_0) = 8\pi, 
\end{equation}
see \Cref{def:Omega},
one may use suitably scaled axisymmetric embedded spheres and tori constructed in \cite[Theorem 1.2. and 1.1]{Scharrer22NLA} to see that for all $c_0<0$ and for all $\sigma\in(0,1/2)$ there exist $\bar A_0,\bar V_0>0$ with $36\pi \bar V_0^2/\bar A_0^3 = \sigma$ such that for all $0<A_0<\bar A_0$, $0<V_0<\bar V_0$ with $36\pi V_0^2/A_0^3 = \sigma$, Condition \eqref{eq:en-thresh} admits axisymmetric surfaces $f_0$. By uniqueness, the evolution equation \eqref{eq:intro:gradient-flow} preserves axial (and reflection) symmetry, and thus we obtain the following.
\begin{corollary}
Let $\Sigma\in\{\mathbb S^2,\mathbb T^2\}$ and suppose that $A_0,V_0>0$ satisfy \eqref{eq:iso}. There exists $\varepsilon^\Sigma(A_0,V_0)>0$ and $\bar A_0^\Sigma,\bar V_0^\Sigma>0$ such that if $c_0\in \R$ and  
\begin{align}
|c_0|<\varepsilon^\Sigma(A_0,V_0) \qquad \text{ or } \qquad c_0<0,\; A_0<\bar{A}_0^\Sigma, \; V_0<\bar{V}_0^\Sigma, \; 36\pi V_0^2/A_0^3< 1/2,
\end{align}
then there exists an axisymmetric constrained $\mathcal{H}_{c_0}$-Helfrich immersion $f_\infty \colon \Sigma\to\R^3$, which is symmetric with respect to a plane orthogonal to the axis of revolution, and satisfies $\mathcal{A}(f_\infty)=A_0$ and $\mathcal{V}(f_\infty)=V_0$.
\end{corollary}

Among axisymmetric surfaces of genus $g\in\{0,1\}$, existence of minimizers for \eqref{intro:HP} has been shown in \cite{choksiveneroni2013}. However, these minimizers do not need to be smooth, and thus our findings provide the first axisymmetric spheres and tori solving \eqref{eq:intro:Helfrich_eq} for $c_0\neq 0$. In the case $c_0=0$, it follows from 
\cite{metsch2023axially} that the minimizers of \cite{choksiveneroni2013} are smooth.

\subsection{Relation to previous works}
For general $c_0\in \R$, global existence of an energy-decreasing evolution for $\mathcal{H}_{c_0}$ that preserves both area and volume has recently been proven in \cite{BKS2024}. Their approach is based on the theory of varifolds, a weak notion of surfaces, and De Giorgi’s \emph{Generalized Minimizing Movements.} While allowing for very general and possible singular shapes, the physical interpretation of these varifold solutions is limited, since it is still unknown whether solutions that are constant in time or their limits as $t\to\infty$ are equilibria of the system. Instead, \Cref{thm:main-result,intro:cor:axisymmetric_tori} provide a dynamic version of the Canham--Helfrich model \eqref{intro:HP} that yields classical smooth solutions, converging to equilibria.

In contrast to the previous work \cite{Lengeler2018} 
on local stability of minimizers under the flow, in \Cref{thm:main-result,intro:cor:axisymmetric_tori} we do not assume the initial datum to be close, in parametrization, to a minimizer. Instead, \eqref{eq:en-thresh} gives an explicit threshold, depending on the parameters, below which we have global existence and convergence. This is a common assumption for Willmore-type gradient flows, see \cite{kuwertschaetzle2004,Jachan,dallacquamullerschatzlespener2020,rupp2023,Jakob2023,rupp2024}, which is also necessary, in general, by Blatt's counterexample \cite{blatt2009}. In \cite{schlierf2024helfrich}, the third author considers the nonprojected analog of \eqref{eq:intro:proj_gradient-flow} where area and volume are not preserved along the evolution. While a positive result is obtained for $c_0>0$, in contrast to \Cref{thm:main-result}, in the physical case $c_0<0$ \cite{DeulingHelfrich} the evolution of the nonprojected flow can shrink to a ``round point'' in finite time, see also \cite{McCoyWheeler2016,Blatt2019} for the case $c_0=0$. Thus, fixing area and volume along the flow prevents such undesirable finite-time singularities and eventually produces constrained critical points of $\H_{c_0}$ as a limit.

In literature, there are various computational approaches to the area- and volume-preserving Helfrich flow \cite{barrettgarckenuernberg2008,barrettgarckenuernberg2016,barrettgarckenuernberg2019,barrettgarckenuernberg2021}. For instance, numerical simulations of $(A_0,V_0,0)$-Helfrich flows can be found in \cite[Figures~14--17]{barrettgarckenuernberg2008} and \cite[Figure~5]{barrettgarckenuernberg2021}, where bi-concave surfaces are numerically obtained as limiting shapes of the flow already for $c_0=0$. 

Sharpness of the energy threshold \eqref{eq:en-thresh} is still an open problem. For the Willmore flow, the existence of singularities was first suggested by numerical experiments in \cite{mayersimonett2002} and only analytically verified years later in \cite{blatt2009}, using the classification of 
axisymmetric Willmore surfaces in \cite{langersinger1984}. Asymptotically, the initial data for the Willmore flow considered in the singular examples in \cite{mayersimonett2002,blatt2009} approach the varifold $V_{\mathrm{lim}}$ made up of two round spheres of radius $1$ touching in a single point, up to scaling. One finds $A_0=\A(V_{\mathrm{lim}})=8\pi$, $V_0=\V(V_{\mathrm{lim}})=\frac83\pi$ and $\sigma=\frac12$, using \eqref{eq:iso}. Thus, for $c_0\in (-\infty,2]$,
\begin{equation}
    \sqrt{\H_{c_0}(V_{\mathrm{lim}})} = \sqrt{8\pi}-\frac{c_0}2\sqrt{8\pi} \geq \min\Bigl\{\sqrt{\Omega(A_0,V_0,c_0)}, \sqrt{\frac{4\pi}{\sigma}} - \frac{c_0}2\sqrt{A_0}\Bigr\}
\end{equation}
with equality for $c_0\in [0,2]$ so that our energy threshold \eqref{eq:en-thresh} is attained 
in this case. Whether the initial data considered in \cite{mayersimonett2002,blatt2009} actually result in singular $(A_0,V_0,c_0)$-Helfrich flows is still an open problem. However, for the nonprojected analog of \eqref{eq:intro:proj_gradient-flow} analytically treated in \cite{schlierf2024helfrich}, there is numerical evidence in this direction. Indeed, numerical experiments in \cite[Section~5.3, Figure~9]{barrettgarckenuernberg2016} suggest that, for $c_0=2$, starting in a suitable spherical ``tube'', one has convergence (in a weak sense) of the nonprojected analog of \eqref{eq:intro:proj_gradient-flow} to $V_{\mathrm{lim}}$.

\subsection{New challenges and strategy of the proofs}
In the singularity analysis of \eqref{eq:intro:gradient-flow}, there are several new difficulties to overcome compared to previous works on (constrained) \emph{Willmore} flows, cf.\ in particular \cite{rupp2024}. First, in general, the classical Li--Yau inequality 
\eqref{eq:intro:LY}
cannot be used to conclude embeddedness from smallness of the Helfrich energy, see \cite[Example 1.2]  {ruppscharrer2023}. 
Fortunately, the perfectly suited Li--Yau-type inequality \eqref{eq:intro:HLY} was proven in \cite{ruppscharrer2023} for exactly this purpose.
Since the limit of the flow in \Cref{thm:main-result} is in particular a smooth constrained critical point for \eqref{intro:HP}, the parameters $A_0,V_0,c_0$ should be in the same range as in the corresponding stationary problem \eqref{intro:eq:energy-condition}, hence \eqref{eq:en-thresh} comprises the assumption
\begin{equation}\label{eq:en-threshold-HLY}
    \mathcal H_{c_0}(f_0) \le \Omega(A_0,V_0,c_0).
\end{equation}
Second, the degeneracy of the constraints in \eqref{eq:intro:gradient-flow} needs to be controlled by proving a uniform lower bound on \eqref{eq:CMC-deficit}. This corresponds to a \emph{quantitative Aleksandrov inequality} (see \cite{Ciraolo2021,julin2024sharpquantitativealexandrovinequality}) which cannot be simply achieved by the triangle inequality in contrast to \cite[Lemma 4.1]{rupp2024}. Third, to obtain a \emph{localized curvature concentration-based criterion} \`a la Kuwert--Schätzle \cite{kuwertschaetzle2002}, appropriate \emph{a priori estimates} that reduce the order of the scaling-critical nonlocal Lagrange multipliers in \eqref{eq:intro:gradient-flow} are needed. 
For this purpose, the condition \eqref{eq:en-thresh} ensures
\begin{equation}\label{eq:en-threshold_lagrange}
    \sqrt{\mathcal H_{c_0}(f_0)}\le \sqrt{\frac{4\pi}{\sigma}} - \frac{c_0}2\sqrt{A_0},
\end{equation}
simplifying
to $\mathcal{W}(f_0)\leq 4\pi/\sigma$ if $c_0=0$ which is in accordance with \cite{rupp2024}. Fourth, we would like to use the removability result in \cite{kuwertschaetzle2004} to apply \cite{bryant1984} in the analysis of noncompact parabolic blow-up limits. However, this requires control over quantities appearing in the monotonicity formula for the Willmore energy \cite{simon1993,kuwertschaetzle2004} that can generically not be controlled by the Helfrich energy. 
Last, it is crucial to preserve orientation during all reparametrization processes to ensure the invariance of both the energy and the flow.

To prove \Cref{thm:main-result}, we first follow the methods of Kuwert and Schätzle \cite{kuwertschaetzle2001,kuwertschaetzle2002} and study a localized version of the energy in which the Lagrange multipliers appear, see \Cref{sec:energy}. Exploiting the {scaling behavior} of the involved geometric quantities, we obtain a linear system for $\lambda_1,\lambda_2$, see \Cref{lem:lag}. Under the assumption \eqref{eq:en-thresh}, this allows us to reduce the order of the Lagrange multipliers to obtain a concentration-based lifespan bound similar to \cite{kuwertschaetzle2002} in \Cref{prop:life}.
We then follow \cite{kuwertschaetzle2001,kuwertschaetzle2004} to construct a blow-up, which we can identify as a (possibly noncompact) and nonplanar Helfrich immersion. Crucially, we find that whenever this blow-up limit is noncompact, it has to be a Willmore immersion, see \Cref{prop:concentration-limit}, whereas if it is compact, we may argue as in \cite{chillfasangovaschaetzle2009} to conclude convergence employing an appropriate {\L}ojasiewicz inequality \cite{Rupp2020} in \Cref{sec:loja}. In \Cref{sec:conv_sphere}, we exclude the possibility of noncompact blowups under the assumption \eqref{eq:en-thresh}. This follows from a precise examination of the scaling behavior of the monotonicity formula for the Helfrich functional from \cite{ruppscharrer2023}, which enables us to apply the removability result from \cite{kuwertschaetzle2004} and then Bryant's classification \cite{bryant1984}. Lastly, in \Cref{sec:tori}, we prove \Cref{intro:cor:axisymmetric_tori}. To that end, we use the relation with elastic curves in the hyperbolic plane as in \cite{dallacquamullerschatzlespener2020} and the observation that unbounded hyperbolic length necessarily produces a self-intersection on the axis of revolution \cite{schlierf2024DirichletWillmore} which is excluded by the Li--Yau inequality 
\eqref{eq:intro:HLY}

\section{Geometric preliminaries}\label{sec:prelim}

\subsection{Basic definitions} \label{sec:def}

Throughout this article, $\Sigma$ denotes a closed, oriented and connected surface. For an immersion $f\colon\Sigma\to\R^3$ of $\Sigma$, $g\vcentcolon =f^*\langle\cdot,\cdot\rangle$ denotes the pull-back metric on $\Sigma$, i.e., $g_{ij}\vcentcolon=\langle\partial_i f,\partial_j f\rangle$ in local coordinates, with induced measure $\mu$. Further, $\nu\colon\Sigma\to\S^2$ is the unique smooth unit normal field induced by the orientation of $\Sigma$. In local coordinates with respect to the orientation of $\Sigma$, one has
\begin{equation}\label{eq:def-nu}
  \nu\vcentcolon=\frac{\partial_1 f\times\partial_2 f}{|\partial_1 f\times \partial_2 f|}.
\end{equation}
The (scalar) second fundamental form of $f$ is the $(2,0)$-tensor determined by $A_{ij}\vcentcolon=\langle\partial^2_{ij}f,\nu\rangle$ in local coordinates. Its trace is the mean curvature $H\vcentcolon=g^{ij}A_{ij}$ where $(g^{ij}) \vcentcolon= (g_{ij})^{-1}$. The trace-free second fundamental form is given by $A_{ij}^0\vcentcolon =A_{ij}-\frac12 H g_{ij}$. Its squared $L^2(\dd\mu)$-norm is denoted by $\W_0(f)\vcentcolon=\int_{\Sigma}|A^0|^2\dd \mu$. Using these definitions, one finds the relation
\begin{align}
  |A|^2&=|A^0|^2+\frac12 H^2.\label{eq:a-a0-h}
\end{align}
For any $c_0\in\R$, define the Helfrich energy with spontaneous curvature $c_0$
\begin{equation}
  \H_{c_0}(f)\vcentcolon=\frac{1}{4}\int_{\Sigma}(H-c_0)^2\dd\mu
\end{equation}
where in the special case $c_0=0$, $\W(f)=\H_{0}(f)$ is also called Willmore energy. Note that, as a consequence of the Gauss-Bonnet theorem, using that the Gauss curvature $K$ can be written as $K=\frac14H^2 - \frac12|A^0|^2$, one has
\begin{equation}\label{eq:W_0=W+g}
  \W_0(f)= \int_{\Sigma}|A^0|^2\dd\mu = 2\W(f)-8\pi(1-\genus\Sigma).
\end{equation}
One further defines
\begin{equation}
    \A(f)=\mu(\Sigma),\quad \V(f)=-\frac13\int_{\Sigma}\langle f,\nu\rangle\dd\mu,\quad \mathcal{I}(f)=36\pi \frac{\V(f)^2}{\A(f)^3},
\end{equation}
the area, signed volume, and isoperimetric ratio of $f$, respectively. Throughout this article, the orientation of $\Sigma$ is chosen such that $\V(f)\geq 0$. In particular, if $f$ is an embedding, $\nu$ is inward pointing so that round spheres of radius $r>0$ have constant mean curvature $H\equiv 2/r>0$.

\begin{lemma}\label{lem:will-vs-helf}
  Let $c_0\in\R$ and consider an immersion $f\colon\Sigma\to\R^3$. Then
  \begin{equation}\label{eq:will-vs-helf}
    \W(f)\leq \Bigl(\sqrt{\H_{c_0}(f)}+\frac12\sqrt{c_0^2\A(f)}\Bigr)^2.
  \end{equation}
\end{lemma}
\begin{proof}
  Without loss of generality, suppose $c_0\neq 0$. By a direct computation, one finds for any $\varepsilon\in(0,1)$
  \begin{align}
    \W(f) &= \H_{c_0}(f) + \frac12 c_0\int_{\Sigma} H\dd\mu - \frac14c_0^2\A(f) \leq \H_{c_0}(f) + \sqrt{\W(f)}\sqrt{c_0^2\A(f)} - \frac14c_0^2\A(f)\\
    &\leq \H_{c_0}(f) + \varepsilon\W(f) + \frac14\frac{1-\varepsilon}{\varepsilon} c_0^2\A(f),
  \end{align}
   by Young's inequality. This yields
  \begin{equation}
    \W(f) \leq \frac{1}{1-\varepsilon} \H_{c_0}(f) + \frac{1}{4\varepsilon} c_0^2\A(f).
  \end{equation}
  Minimizing with respect to $\varepsilon\in(0,1)$ yields the claim.
\end{proof}

\subsection{Evolution equations}
First, recall the following well-known evolutions of the relevant geometric quantities.
\begin{lemma} For a family of immersions $f\colon[0,T)\times \Sigma \to\R^3$ with normal speed $\partial_t f=\vcentcolon \xi\nu$, one has in the coordinates of a local orthonormal frame ${e_1,e_2}$ 
  \begin{align}
    \partial_t \dd\mu &= -H\xi\dd\mu\label{eq:ev-mu}\\
    \partial_t H &= \Delta\xi + |A|^2\xi=\Delta\xi + |A^0|^2\xi + \frac{1}{2}H^2\xi\label{eq:ev-H}\\
  \end{align}
\end{lemma}

As in \cite[Proposition 2.4]{rupp2024}, these evolutions immediately yield

\begin{proposition}\label{prop:variations}
  Let $f\colon\Sigma\to\R^3$ be an immersion and $\varphi\colon\Sigma\to\R^3$. One has the following first variation identities.
  \begin{align}
    \W_0'(f)(\varphi)&=\int_{\Sigma} \langle [\Delta H+|A^0|^2H]\nu,\varphi\rangle\dd\mu =\vcentcolon \langle \nabla\W_0(f),\varphi\rangle_{L^2(\dd\mu)},\\
    \A'(f)(\varphi)&=-\int_{\Sigma} \langle H\nu,\varphi\rangle\dd\mu =\vcentcolon \langle \nabla\A(f),\varphi\rangle_{L^2(\dd\mu)},\\
    \V'(f)(\varphi)&=-\int_{\Sigma} \langle \nu,\varphi\rangle\dd\mu =\vcentcolon \langle \nabla\V(f),\varphi\rangle_{L^2(\dd\mu)}
  \end{align}
  and finally $\H_{c_0}'(f)(\varphi) = \langle \nabla\H_{c_0}(f),\varphi\rangle_{L^2(\dd\mu)}$ with
  \begin{align}
    2\nabla \H_{c_0}(f)&=\nabla\W_0(f)-c_0|A^0|^2\nu + \frac12c_0H(H-c_0)\nu\\
    &= [\Delta H + |A^0|^2(H-c_0) + \frac12 c_0 H(H-c_0)]\nu. 
  \end{align}
\end{proposition}

\subsection{The area and volume preserving Helfrich flow}\label{sec:preserving_Helfrich_flow}

Given $A_0,V_0>0$ with \eqref{eq:iso} and $c_0\in\R$, a family of immersions $f\in C^\infty([0,T)\times \Sigma;\R^3)$ satisfying $A_0=\A(f(0))$, $V_0=\V(f(0))$ and

\begin{align}
  \partial_t f &= -\Bigl[ \Delta H + |A^0|^2H +  c_0(\frac12H^2-|A^0|^2) - (\lambda_1+\frac12c_0^2) H - \lambda_2 \Bigr]\nu \\
  &= -2 \nabla \mathcal{H}_{c_0}(f) + \lambda_1(f) H\nu + \lambda_2(f) \nu \label{eq:flow-eq}
\end{align}
where 
\begin{equation}\label{eq:def-l1}
  \lambda_1\vcentcolon = \lambda_1(f)=2\frac{\A(f)\int_{\Sigma}\langle\nabla\H_{c_0}(f),H\nu\rangle\dd\mu-\int_{\Sigma}H\dd\mu\int_{\Sigma}\langle\nabla\H_{c_0}(f),\nu\rangle\dd\mu}{4\W(f)\A(f)-(\int_{\Sigma}H\dd\mu)^2}
\end{equation}
and
\begin{equation}\label{eq:def-l2}
  \lambda_2\vcentcolon =\lambda_2(f)=2\frac{-\int_{\Sigma}H\dd\mu\int_{\Sigma}\langle\nabla\H_{c_0}(f),H\nu\rangle\dd\mu+4\W(f)\int_{\Sigma}\langle\nabla\H_{c_0}(f),\nu\rangle\dd\mu}{4\W(f)\A(f)-(\int_{\Sigma}H\dd\mu)^2}
\end{equation}
is referred to as an $(A_0,V_0,c_0)$-Helfrich flow starting in $f_0\vcentcolon=f(0)$.

\begin{remark}\label{rem:energy-decay}
  By direct computation, if $f\colon[0,T)\times\Sigma\to\R^3$ is an $(A_0,V_0,c_0)$-Helfrich flow, then $\A(f(t))=A_0$ and $\V(f(t))=V_0$ for all $t\in[0,T)$. Moreover, 
  \begin{equation}\label{eq:en-dec}
    2\frac{\dd}{\dd t}\H_{c_0}(f(t)) =  2 \int_{\Sigma} \langle \nabla\H_{c_0}(f),\partial_tf\rangle\dd \mu  =-\int_{\Sigma} |\partial_tf|^2\dd \mu \leq 0.
  \end{equation}
  Particularly, $t\mapsto\H_{c_0}(f(t))$ is nonincreasing and thus $\lim_{t\nearrow T}\H_{c_0}(f(t))$ exists. Moreover, for a \emph{positive}, that is, orientation preserving diffeomorphism $\phi\colon\Sigma\to\Sigma$, $\tilde{f}\colon[0,T)\times\Sigma\to\R^3$ with $\tilde{f}(t,x)\vcentcolon=f(t,\phi(x))$ again is an $(A_0,V_0,c_0)$-Helfrich flow.
\end{remark}

\begin{remark}
  Writing $\overline{H}=\fint_{\Sigma} H\dd\mu$, we observe that the denominator in \eqref{eq:def-l1} and \eqref{eq:def-l2} satisfies
  \begin{align}
    4\W(f)&{ }\A(f)-\Bigl(\int_{\Sigma}H\dd\mu\Bigr)^2 = \frac12 \int_{\Sigma}\int_{\Sigma} \bigl(H(x)-H(y)\bigr)^2\dd\mu(x)\dd\mu(y) \\
    &= \A(f)\Big(\int_{\Sigma}H^2\dd\mu-\A(f)\overline{H}^2\Big) = \A(f)\int_{\Sigma}(H-\overline{H})^2\dd\mu
  \end{align}
  so that,  with $\overline{\mathcal{H}}(f)$ as in \eqref{eq:CMC-deficit},
  \begin{equation}\label{eq:denom-cmc-estimate}
    4\W(f)\A(f)-\Bigl(\int_{\Sigma}H\dd\mu\Bigr)^2 = 4\A(f)\overline{\H}(f).
  \end{equation}
  By the famous Aleksandrov theorem \cite{Aleksandrov1962}, this vanishes for embedded surfaces only in case of the round sphere. This however is impossible if $36\pi \V(f)^2<\A(f)^3$ by the isoperimetric theorem.
\end{remark}

We now define the constant $\Omega$ that appears in \Cref{thm:main-result}.

\begin{definition}\label{def:Omega}
    For all $c_0\in\R$ and $A_0,V_0>0$ define 
    \begin{equation}
        \Omega(A_0,V_0,c_0)\vcentcolon=
        \begin{cases}
            4\pi+\sqrt{(4\pi)^2+\frac{|c_0|V_0}{2C^2A_0}}&\text{if $c_0<0$} \\
            \max\{(\sqrt{8\pi}-\frac12\sqrt{c_0^2A_0})_+^2,8\pi-6c_0(4\pi^2V_0)^{\frac{1}3}\}&\text{if $c_0\ge0$}
        \end{cases}
    \end{equation}
    where $C$ is the best constant in \cite[Theorem 1.2]{Scharrer23}, and $(x)_+\vcentcolon=\max\{0,x\}$ denotes the positive part of $x\in\R$.
\end{definition}

The constant $\Omega(A_0,V_0,c_0)$ is chosen to control the right hand side in the Li--Yau-type inequality \eqref{eq:intro:HLY}.

\begin{lemma}\label{lem:Omega_motivation}
Let $\delta>0$, $f\colon \Sigma\to\R^3$ be an embedding, and $A_0\vcentcolon=\mathcal{A}(f)$, $V_0\vcentcolon=\mathcal{V}(f)$ satisfy \eqref{eq:iso}. Suppose that $\mathcal{H}_{c_0}(f)\le\Omega(A_0,V_0,c_0)-\delta$. If $c_0\geq 0$ and $\Omega(A_0,V_0,c_0)= (\sqrt{8\pi}-\frac{1}{2}\sqrt{c_0^2A_0})_+^2$, then
\begin{equation}\label{eq:Omega_motivation2}
\frac{1}{4\pi}\mathcal W(f) \le 2-\frac{\delta}{4\pi}.
\end{equation}
Otherwise
\begin{align}\label{eq:Omega_motivation0}
\frac{1}{4\pi} \mathcal{H}_{c_0}(f) - \frac{c_0}{2\pi}\int_\Sigma \frac{\langle f-x_0,\nu\rangle}{|f-x_0|^2}\dd\mu\le2-\frac{\delta}{4\pi}.
\end{align}
\end{lemma}

\begin{proof}
Suppose $c_0\geq 0$ and $\Omega(A_0,V_0,c_0)= (\sqrt{8\pi}-\frac{1}{2}\sqrt{c_0^2A_0})_+^2$. Noting that $\Omega(A_0,V_0,c_0)>\mathcal H_{c_0}(f)\geq 0$,
\Cref{lem:will-vs-helf} implies \eqref{eq:Omega_motivation2}.

If $c_0\ge0$ and $\Omega(A_0,V_0,c_0)= 8\pi - 6 c_0 (4\pi^2 V_0)^{\frac{1}{3}}$, then \cite[Remark 6.11(ii)]{ruppscharrer2023} implies \eqref{eq:Omega_motivation0}. 

Lastly, for $c_0\leq 0$, we use the divergence theorem (cf.\ \cite[(6.18)]{ruppscharrer2023}) combined with the diameter estimate \cite[Theorem 1.2]{Scharrer23}, to infer
\begin{equation}
    -\int_\Sigma \frac{\langle f-x_0,\nu\rangle}{|f-x_0|^2}\dd\mu \ge \frac{V_0}{4C^2 A_0\mathcal H_{c_0}(f)}.
\end{equation}
Hence, since $c_0\leq 0$, we estimate
\begin{equation}\label{eq:cvolume}
\frac{1}{4\pi} \mathcal{H}_{c_0}(f) - \frac{c_0}{2\pi}\int_\Sigma \frac{\langle f-x_0,\nu\rangle}{|f-x_0|^2}\dd\mu\le \frac{\Omega-\delta}{4\pi} + \frac{c_0}{2\pi}\frac{V_0}{4C^2 A_0 \Omega}.
\end{equation}
A straightforward computation shows that with $\Omega=\Omega(A_0,V_0,c_0)$ as in \Cref{def:Omega} for $c_0<0$, we have
\begin{align}
\frac{\Omega}{4\pi} + \frac{c_0}{2\pi}\frac{V_0}{4C^2 A_0 \Omega} = 2,
\end{align}
and \eqref{eq:Omega_motivation0} follows.
\end{proof}

\begin{proposition}\label{prop:scaling_Omega}
    Let $r>0$. For $\Omega(A_0,V_0,c_0)$ as in \Cref{def:Omega}, we have 
    \begin{equation}
        \Omega\Big(\frac{A_0}{r^2},\frac{V_0}{r^3},rc_0\Big) = \Omega(A_0,V_0,c_0).
    \end{equation}
    Moreover, for any smoothly 
    immersed
    oriented surface $f\colon\Sigma\to\R^3$, there holds
    \begin{equation}\label{eq:scaling_Helfrich}
        \H_{rc_0}(r^{-1}f)=\H_{c_0}(f).
    \end{equation}
\end{proposition}

\begin{proposition}\label{prp:Omega-and-LiYau}
    Let $f_0\colon\Sigma\to\R^3$ be a smooth embedding with $A_0=\mathcal A(f_0)$ and $V_0=\mathcal V(f_0)$ satisfying
    \begin{equation}
        \mathcal H_{c_0}(f_0)<\Omega(A_0,V_0,c_0)
    \end{equation}
    for some $c_0\in\R$. Let $f\colon[0,T)\times\Sigma\to\R^3$ be an $(A_0,V_0,c_0)$-Helfrich flow starting in $f_0$. Then $f(t,\cdot)\colon \Sigma\to\R^3$ is an embedding for all $t\in[0,T)$.
\end{proposition}
\begin{proof}
    Let 
    \begin{equation}
        I\vcentcolon=\{t\in[0,T)\mid f(t,\cdot)\text{ is an embedding}\}.
    \end{equation}
    A simple proof by contradiction using \cite[Lemma 6.19]{AndrewsChowGuentherLangford} shows that $I$ is open. Moreover, $I\neq\emptyset$ since $0\in I$. Thus, it is enough to show that $I$ is closed. To this end, let $t_0 \in (0,T)$ be in the closure of $I$. Let $V_t$ with $t\in I$ be the oriented varifolds induced by $f(t,\cdot)$, cf.\ \cite[Example 2.4]{ruppscharrer2023}. By the Jordan--Brouwer separation theorem, all $V_t$ are volume varifolds in the sense of \cite[Definition 3.1]{Scharrer23}, cf.\ \cite[Lemma 6.6]{ruppscharrer2023}. Combining \eqref{eq:a-a0-h}\eqref{eq:W_0=W+g} with \Cref{lem:will-vs-helf} and \Cref{rem:energy-decay}, we see that the second fundamental forms of $V_t$ remain bounded in $L^2$ uniformly in $t$. Thus, by \cite[Theorem 7.6]{Scharrer23} there exists a varifold limit $V_{\infty}$ for a subsequence $t_k\to t_0$ which is again a volume varifold satisfying
    \begin{equation}
        \mathcal H_{c_0}(V_\infty)<\Omega(A_0,V_0,c_0).
    \end{equation}
    Since $f(t_k,\cdot)\to f(t_0,\cdot)$ in $C^1(\Sigma)$, we see that $V_\infty$ is the varifold induced by $f(t_0,\cdot)$. Now, \cite[Corollary 4.11]{ruppscharrer2023} implies that $f(t_0,\cdot)$ is an embedding. Consequently, $t_0\in I$ and $I$ is closed.
\end{proof}

Using \eqref{eq:denom-cmc-estimate}, we control the denominators in \eqref{eq:def-l1} and \eqref{eq:def-l2} with the following
\begin{proposition}\label{prop:cmc-control}
  Let $c_0\in\R$, $A_0,V_0>0$ satisfy \eqref{eq:iso}, and $\varepsilon>0$. Then, denoting the set of smooth embeddings $f\colon\Sigma\to\R^3$ with $\mathcal{A}(f)=A_0$ and $\mathcal V(f)=V_0$ by $\mathcal S_\Sigma(A_0,V_0)$, there holds 
  \begin{equation}
    \bar \eta_\varepsilon(A_0,V_0,c_0,\genus\Sigma)\vcentcolon=\inf\Bigl\{ \overline{\H}(f) \mid f\in \mathcal S_\Sigma(A_0,V_0),\,\mathcal H_{c_0}(f)\leq \Omega(A_0,V_0,c_0)-\varepsilon \Bigr\}>0.
  \end{equation}
\end{proposition}
\begin{proof}
    For the sake of contradiction, assume the existence of a sequence $f_k\in\mathcal S_\Sigma(A_0,V_0)$ such that
    \begin{equation}\label{eq:prop:denominator:energy_bound}
        \sup_{k\in\N}\mathcal H_{c_0}(f_k)\le \Omega(A_0,V_0,c_0) - \varepsilon
    \end{equation}
    and 
    \begin{equation}\label{eq:prop:denominator:contrad_assump}
        \lim_{k\to\infty}\overline \H(f_k)=0.
    \end{equation}
    The total mean curvature is uniformly bounded:
    \begin{equation}\label{eq:prop:denominator:boundedness_barH}
        A_0|\overline H_k| = \Big|\int_\Sigma H_k\dd\mu_k\Big|\le \sqrt{4A_0\Omega(A_0,V_0,c_0)} + |c_0|A_0.
    \end{equation}
    Thus, by passing to a subsequence, we achieve that the limit
    \begin{equation}
        \lim_{k\to\infty}\overline H_k=\vcentcolon h_0\in\R
    \end{equation}
    exists. Moreover,
    \begin{equation}
        \H_{h_0}(f_k)=\overline \H(f_k) + \frac12\big(\overline H_k-h_0\big)\int_{\Sigma}H_k\dd\mu_k - \frac14\big(\overline H_k^2- h_0^2\big)A_0.
    \end{equation}
    Now, we may apply \cite[Theorem 7.6]{Scharrer23} to deduce the existence of a volume varifold $V$ in $\R^3$ with $\mathcal A(V)=A_0$, $\mathcal V(V)=V_0$, and
    \begin{equation}
        \H_{h_0}(V) = 0,\qquad \H_{c_0}(V) \le \Omega(A_0,V_0,c_0)-\varepsilon.
    \end{equation}
    By \cite[Corollary 4.11]{ruppscharrer2023} and \cite[Theorem 5.8]{RuppScharrer24}, the varifold $V$ is induced by a $W^{2,2}$-conformal Lipschitz embedding $F\colon \Sigma_0\to\R^3$ of a Riemann surface $\Sigma_0$. Hence, the property $\H_{h_0}(F) = 0$ implies that in local conformal coordinates defined on the disk $D\subset \R^2$,
    \begin{equation}\label{eq:cmc-equation}
        \Delta F = h_0\partial_1F\times\partial_2F, 
    \end{equation}
    where $\Delta$ is the flat Laplacian in $\R^2$.
    The right hand side is in $W^{1,2}(D)$ since $F$ is Lipschitz continuous. Bootstrapping Equation \eqref{eq:cmc-equation}, it follows that $F$ is smooth by standard elliptic regularity, cf.\ Theorems 5.6.6 and 6.3.2 in \cite{Evans98AMS}. Thus, 
    Aleksandrov's theorem \cite{Aleksandrov1962} implies that $F$ parametrizes a round sphere, contradicting the assumption $36\pi V_0^2/A_0^3 < 1$.    
\end{proof}

The following scaling behavior of the flow and especially the Lagrange multipliers $\lambda_1$ and $\lambda_2$ is immediate from their respective definitions. Also cf. \cite[Lemma 2.8]{rupp2024}.

\begin{lemma}\label{lem:par-scal}
  Let $c_0\in\R$, consider $A_0,V_0>0$ satisfying \eqref{eq:iso} and let $r>0$. If $f\colon[0,T)\times\Sigma\to\R^3$ is an $(A_0,V_0,c_0)$-Helfrich flow, then its parabolic rescaling $\widetilde{f}\colon[0,T/r^4)\times\Sigma\to\R^3$ with $\widetilde{f}(t,p)=\frac{1}{r}f(r^4t,p)$ satisfies the following.
  \begin{enumerate}[(i)]
    \item $\widetilde{f}$ is an $(\frac{1}{r^2}A_0,\frac{1}{r^3}V_0,rc_0)$-Helfrich flow.
    \item $\widetilde{\lambda_1}(t)=r^2\lambda_1(r^4t)$ and $\widetilde{\lambda_2}(t)=r^3\lambda_2(r^4t)$.
    \item $\int_0^T(\lambda_1(t))^2\dd t =\int_{0}^{T/r^4}(\widetilde{\lambda_1}(t))^2\dd t$ and $\int_0^T|\lambda_2(t)|^{\frac{4}{3}}\dd t =\int_{0}^{T/r^4}|\widetilde{\lambda_2}(t)|^\frac{4}{3}\dd t$.
  \end{enumerate}
\end{lemma}

\section{Localized energy estimates}\label{sec:energy}

This section is devoted to the fundamental estimates which enable a suitable blow-up construction in later sections. To this end, following \cite[Section 3 respectively]{kuwertschaetzle2001,rupp2023,rupp2024}, one localizes the energy decay of $\int |A|^2\dd \mu$ to obtain a life-span theorem and suitable uniform bounds on derivatives of the curvature as long as an energy concentration is controlled.

Similarly as in \cite{kuwertschaetzle2001,rupp2023,rupp2024}, consider $\widetilde{\gamma}\in C_c^{\infty}(\R^3)$ with $|\widetilde{\gamma}|\leq 1$ and define $\gamma=\widetilde{\gamma}\circ f$. Fix some $\Lambda>0$ 
such that $\|D\widetilde{\gamma}\|_{\infty}\leq\Lambda$ and $\|D^2\widetilde{\gamma}\|_{\infty}\leq\Lambda^2$. As on \cite[p.~7]{rupp2023}, one estimates
\begin{equation}\label{eq:gamma}
  |\nabla\gamma|\leq \Lambda\quad\text{and}\quad |\nabla^2\gamma|\leq C(\Lambda^2+|A|\Lambda).
\end{equation}

In the following, $C\in (0,\infty)$ is a universal constant that may change from line to line.

\begin{lemma}\label{lem:loc-2}
    Let $A_0,V_0>0$ satisfy \eqref{eq:iso}, $c_0\in\R$ and $f\colon[0,T)\times\Sigma\to\R^3$ be an $(A_0,V_0,c_0)$-Helfrich flow. With $\tilde{\gamma}$, $\gamma$ as in \eqref{eq:gamma}, we have
    \begin{align}
    \frac{\dd}{\dd t}&\int_{\Sigma} |A|^2\gamma^4\dd\mu + \frac32\int_{\Sigma}|\nabla\W_0(f)|^2\gamma^4\dd\mu\\
    &\leq C(|\lambda_1|+c_0^2)\int_{\Sigma}\bigl(|\nabla^2H||A|+|A|^4\bigr)\gamma^4+\Lambda|A|^3\gamma^3\dd\mu+C|\lambda_2|\int_{\Sigma}|A|^{3}\gamma^4+\Lambda|A|^2\gamma^{3}\dd\mu\\
    &\quad + C\Lambda^4\int_{\{\gamma>0\}}|A|^2\dd\mu+C\Lambda^2\int_{\Sigma}|A|^4\gamma^2\dd\mu + C c_0^2\int_{\Sigma}|A|^4\gamma^4\dd\mu.
  \end{align}
\end{lemma}

The proof of \Cref{lem:loc-2} can be found in \Cref{app:local}. The experienced reader can obtain it from
\cite[Lemma 3.2]{rupp2024}, formally replacing $3 \lambda /\mathcal{A}$ with $\lambda_1+c_0^2/2$ and $-2\lambda/\mathcal{V}$ with $\lambda_2$, in combination with \cite[Lemma 4.2]{schlierf2024helfrich} for the additional terms linear in $c_0$.

\begin{proposition}\label{prop:en-est}
  There exist universal constants $\varepsilon_0,\delta_0,C\in(0,\infty)$ such that, for any $f$ and $\gamma$ as in \Cref{lem:loc-2}, if
  \begin{equation}
    \int_{\{\gamma>0\}} |A|^2\dd\mu<\varepsilon_0\quad\text{at some time $t\in[0,T)$},
  \end{equation}
  then at time $t$ one has
  \begin{align}
    \frac{\dd}{\dd t}&\int_{\Sigma}|A|^2\gamma^4\dd\mu + \delta_0 \int_{\Sigma} \bigl( |\nabla^2 A|^2+|A|^2|\nabla A|^2+|A|^6 \bigr)\gamma^4\dd\mu\\
    &\leq C\Lambda^4\int_{\{\gamma>0\}}|A|^2\dd\mu + C(\lambda_1^2+|\lambda_2|^{\frac43}+c_0^4)\int_{\Sigma}|A|^2\gamma^4\dd\mu.
  \end{align}
\end{proposition}
\begin{proof}
  By \cite[Proposition 3.2]{rupp2023}, see also \cite[Proposition 2.6 and Lemma 4.2]{kuwertschaetzle2002},
  at time $t$, the following interpolation inequality holds.
  \begin{equation}
    \int_{\Sigma} \bigl( |\nabla^2 A|^2+|A|^2|\nabla A|^2+|A|^6 \bigr)\gamma^4\dd\mu\leq C\int_{\Sigma}|\nabla\W_0(f)|^2\gamma^4\dd\mu+C\Lambda^4\int_{\{\gamma>0\}}|A|^2\dd\mu.
  \end{equation}
  Using \Cref{lem:loc-2}, there thus exists $\delta_0\in(0,\infty)$ with
  \begin{align}
    \frac{\dd}{\dd t}&\int_{\Sigma} |A|^2\gamma^4\dd\mu + 2\delta_0\int_{\Sigma} \bigl( |\nabla^2 A|^2+|A|^2|\nabla A|^2+|A|^6 \bigr)\gamma^4\dd\mu\\
    &\leq C(|\lambda_1|+c_0^2)\int_{\Sigma}\bigl(|\nabla^2H||A|+|A|^4\bigr)\gamma^4+\Lambda|A|^3\gamma^3\dd\mu+C|\lambda_2|\int_{\Sigma}|A|^{3}\gamma^4+\Lambda|A|^2\gamma^{3}\dd\mu\\
    &\quad + C\Lambda^4\int_{\{\gamma>0\}}|A|^2\dd\mu+C\Lambda^2\int_{\Sigma}|A|^4\gamma^2\dd\mu + C c_0^2\int_{\Sigma}|A|^4\gamma^4\dd\mu.
  \end{align}
  Proceeding as in \cite[Proposition 3.3]{rupp2024}, one finds that
  \begin{align}
    \frac{\dd}{\dd t}&\int_{\Sigma} |A|^2\gamma^4\dd\mu + \frac32\delta_0\int_{\Sigma} \bigl( |\nabla^2 A|^2+|A|^2|\nabla A|^2+|A|^6 \bigr)\gamma^4\dd\mu\\
    &\leq C\Lambda^4\int_{\{\gamma>0\}}|A|^2\dd\mu + C(\lambda_1^2+|\lambda_2|^{\frac43}+c_0^4)\int_{\Sigma}|A|^2\gamma^4\dd\mu + C c_0^2\int_{\Sigma}|A|^4\gamma^4\dd\mu.
  \end{align}
  The claim follows noting that $c_0^2\int_{\Sigma}|A|^4\gamma^4\dd\mu\leq \varepsilon \int_{\Sigma}|A|^6\gamma^4\dd\mu + C(\varepsilon)c_0^4 \int_{\Sigma}|A|^2\gamma^4\dd\mu$.
\end{proof}

Still following the arguments of \cite{kuwertschaetzle2001,rupp2023,rupp2024}, we introduce the following notation.

\begin{definition}
  Consider a family of immersions $f\colon[0,T)\times\Sigma\to\R^3$. For $t\in[0,T)$ and $r>0$, the \emph{curvature concentration function} is defined as 
  \begin{equation}
    \scurv(t,r)=\sup_{x\in\R^3} \int_{B_r(x)} |A|^2\dd\mu.
  \end{equation} 
\end{definition}
Following the notation of \cite{kuwertschaetzle2002}, the integrals above have to be understood over the preimage of the open ball $B_r(x)\subset \R^3$ under $f(t)$, i.e.,  the domain of integration is $(f(t))^{-1}(B_r(x))\subset \Sigma$.
With this notation, one obtains the following corollary as an integrated version of \Cref{prop:en-est}.

\begin{corollary}\label{cor:en-est}
  Let $\varepsilon_0\in(0,\infty)$, $\delta_0\in(0,\infty)$ as in \Cref{prop:en-est}. There exists a uniform constant $C>0$ such that, for $A_0,V_0>0$ satisfying \eqref{eq:iso}, $c_0\in\R$ and an $(A_0,V_0,c_0)$-Helfrich flow $f\colon[0,T)\times\Sigma\to\R^3$, one has the following.

  If there exists $\rho>0$ such that
  \begin{equation}
    \scurv(t,\rho)<\varepsilon_0\quad\text{for all $t\in[0,T)$},
  \end{equation}
  then, for all $x\in\R^3$ and $0\leq t<T$,
  \begin{align}
    &\int_{B_{\frac{\rho}{2}}(x)}|A|^2\dd\mu\Big|_t + \delta_0\int_{0}^{t}\int_{B_{\frac{\rho}{2}}(x)}\bigl( |\nabla^2A|^2+|A|^2|\nabla A|^2+|A|^6 \bigr)\dd\mu\dd\tau\\
    &\leq \int_{B_{\rho}(x)}|A|^2\dd\mu\Big|_{t=0} + \frac{C}{\rho^4} \int_{0}^{t}\int_{B_{\rho}(x)}|A|^2\dd\mu\dd\tau + C\int_{0}^{t}\bigl(\lambda_1^2+|\lambda_2|^{\frac43}+c_0^4\bigr)\int_{B_{\rho}(x)}|A|^2\dd\mu\dd\tau.
  \end{align}
\end{corollary}
\begin{proof}
  Choose $\widetilde{\gamma}\in C_c^{\infty}(\R^3)$ with $\chi_{B_{\frac{\rho}{2}}(x)}\leq \widetilde{\gamma}\leq \chi_{B_{\rho}(x)}$, $\|D\widetilde{\gamma}\|_{\infty}\leq C/\rho$ and $\|D^2\widetilde{\gamma}\|_{\infty}\leq C/\rho^2$. The claim follows by simply integrating \Cref{prop:en-est} in time with $\Lambda=C/\rho$.
\end{proof}

\section{A priori control of the Lagrange multipliers}\label{sec:lagr}

\begin{lemma}\label{lem:lag}
  Consider $A_0,V_0>0$ satisfying $\sigma=36\pi\frac{V_0^2}{A_0^3}\in(0,1)$ and $c_0\in\R$. If $f\colon[0,T)\times\Sigma\to\R^3$ is an $(A_0,V_0,c_0)$-Helfrich flow with 
  \begin{equation}\label{eq:4pi-durch-sigma}
      \sqrt{\H_{c_0}(f_0)}+\frac{c_0}{2}\sqrt{A_0}\leq\sqrt{\frac{4\pi}{\sigma}}(1-\eta)
  \end{equation}
  for some $\eta\in(0,1)$, we find for $j=1,2$
  \begin{align}
    |\lambda_j| &\leq C(\sigma,\eta,\genus\Sigma)\frac{(1+\sqrt{c_0^2A_0})^{j-1}}{A_0^{j/2}} \\
    &\quad \cdot \Bigl(\Bigl[\int_{\Sigma}|\partial_tf|\dd\mu + |c_0|\Bigl]\big(1+\sqrt{c_0^2A_0}\big) + c_0^2\int_{\Sigma}|A|\dd\mu + \int_{\Sigma}|A|^3\dd\mu\Bigr).
  \end{align}
\end{lemma}
\begin{proof}
  The assertion is proved by testing the evolution equation \eqref{eq:flow-eq} in two geometrically natural ways which yields a linear system for $\lambda_1$ and $\lambda_2$.

  \textbf{Affine variation:} Fix any $t\in[0,T)$ and $x_t\in f(t,\Sigma)$. Using \eqref{eq:flow-eq}, \eqref{eq:scaling_Helfrich} and the scaling properties of area and volume,
  \begin{align}
      -\int_{\Sigma} \langle\partial_tf,f-x_t\rangle\dd\mu &= \left.\frac{\dd}{\dd \rho}\right\vert_{\rho=1} \left(2\H_{c_0}(\rho (f-x_t))+\lambda_1(f)\A(\rho(f-x_t))+\lambda_2(f)\V(\rho(f-x_t))\right)\\
      &= 2 \left.\frac{\dd}{\dd \rho}\right\vert_{\rho=1}\H_{\rho c_0}(f) + 2\lambda_1(f)\A(f) + 3\lambda_2(f)\V(f) \\
      &= \left.\frac{\dd}{\dd \rho}\right\vert_{\rho=1} \Big(2\W(f)-\rho c_0 \int_{\Sigma}H\dd\mu+\frac12c_0^2\rho^2\A(f)\Big) + 2\lambda_1(f)\A(f) + 3\lambda_2(f)\V(f)\\
      &=2\lambda_1(f)\A(f) + 3\lambda_2(f)\V(f) - c_0\int_{\Sigma}H\dd\mu +c_0^2\A(f).\label{eq:lem-lag-1}
  \end{align}
  Note that one can estimate $|f-x_t|\leq \mathrm{diam}(f)$.

  \textbf{Variation of the volume:} Multiplying \eqref{eq:flow-eq} by $\nu$ yields at time $t$
  \begin{align}
    0&= \partial_t\V(f) = -\int_{\Sigma}\langle\partial_t f,\nu\rangle\dd\mu \\
    &= \int_{\Sigma}|A^0|^2H\dd\mu-c_0\int_{\Sigma}(|A^0|^2-\frac12H(H-c_0))\dd\mu - \lambda_1\int_{\Sigma}H\dd\mu - \lambda_2\A(f)\\
    &= \int_{\Sigma}|A^0|^2H\dd\mu+2c_0\int_{\Sigma}K\dd\mu  -(\lambda_1+\frac12c_0^2)\int_{\Sigma}H\dd\mu - \lambda_2\A(f)\\
    &=\int_{\Sigma}|A^0|^2H\dd\mu+c_08\pi(1-\genus \Sigma) -(\lambda_1+\frac12c_0^2)\int_{\Sigma}H\dd\mu - \lambda_2\A(f),\label{eq:lem-lag-2}
  \end{align}
  where we used the Gauss--Bonnet theorem.
  That is, we obtain a linear system for $\lambda_1$ and $\lambda_2$ with coefficient matrix $M$ given by
  \begin{equation}
    M=\begin{pmatrix}
      2A_0&3V_0\\
      \int_{\Sigma}H\dd\mu&A_0
    \end{pmatrix}.
  \end{equation}
  First, due to \eqref{eq:4pi-durch-sigma},
  \begin{equation}\label{eq:lagr-a-priori-control-int-H}
    \int_{\Sigma} H\dd\mu \leq 2\sqrt{\H_{c_0}(f)A_0} + c_0A_0 \leq 2\sqrt{\frac{4\pi}{\sigma}}\sqrt{A_0}(1-\eta).
  \end{equation}
  Thus, also using $3V_0=\sqrt{\frac{\sigma}{4\pi}}A_0^{\frac32}$, one finds
  \begin{align}
    \det M &= 2A_0^2-3V_0\int_{\Sigma}H\dd\mu = 2A_0^2\bigl(1-\sqrt{\frac{\sigma}{4\pi}}\frac{1}{2\sqrt{A_0}}\int_{\Sigma}H\dd\mu\bigr)\geq 2A_0^2\eta>0.\label{eq:lem-lag-3}
  \end{align}
  Therefore, the solution of the linear system \eqref{eq:lem-lag-1}, \eqref{eq:lem-lag-2} is 
  \begin{equation}
      \begin{pmatrix}
          \lambda_1 \\ \lambda_2
      \end{pmatrix}
      = \frac{A_0}{\det M}
      \begin{pmatrix}
          1 & - \sqrt{\frac{\sigma}{4\pi}}\sqrt{A_0} \\
          -\fint_\Sigma H\dd\mu & 2
      \end{pmatrix}
      \begin{pmatrix}
          -\int_\Sigma\langle\partial_tf,f-x_t\rangle\dd\mu + c_0\int_\Sigma (H-c_0)\dd\mu \\
          \int_\Sigma|A^0|^2H\dd\mu +c_08\pi(1-\genus \Sigma) -\frac{c_0^2}2\int_\Sigma H\dd\mu
      \end{pmatrix}.
  \end{equation}
  Next, we notice that if $c_0\ge0$, then Simon's diameter bound \cite[Lemma 1.1]{simon1993} and the inequalities \eqref{eq:will-vs-helf}, \eqref{eq:4pi-durch-sigma} imply   
  \begin{equation}\label{eq:diam-bd-c0geq0}
    |f-x_t|\leq \mathrm{diam}(f) \leq C\sqrt{A_0\W(f)} \leq C(\sigma,\eta) \sqrt{A_0}.
  \end{equation}
  On the other hand, if $c_0<0$, then \cite[Theorem 1.2]{Scharrer23} and the inequality \eqref{eq:4pi-durch-sigma} imply 
  \begin{equation}\label{eq:diam-bd-c0leq0}
    |f-x_t|\leq \mathrm{diam}(f) \leq C\sqrt{A_0\mathcal H_{c_0}(f)} \leq C(\sigma,\eta)(1+\sqrt{c_0^2A_0}) \sqrt{A_0}.
  \end{equation}  
  For estimating $|\lambda_2|$, we especially use 
  \begin{equation}
    \Big|\fint_\Sigma H\dd\mu\Big|\leq \frac{1}{{A_0}}\Big(2\sqrt{\H_{c_0}(f)A_0}+|c_0|A_0\Big)\leq \frac{C(\sigma,\eta)}{\sqrt{A_0}}\Big(1+\sqrt{c_0^2A_0}\Big),
  \end{equation} by Cauchy-Schwarz and \eqref{eq:4pi-durch-sigma}. Moreover, again by Cauchy-Schwarz and \eqref{eq:4pi-durch-sigma}, 
  \begin{equation}
    \Big|c_0\int_\Sigma (H-c_0)\dd\mu\Big| \leq 2\sqrt{\H_{c_0}(f)}\sqrt{A_0} |c_0|\leq C(\sigma,\eta)\big(1+\sqrt{c_0^2A_0}\big)\sqrt{A_0}|c_0|.
  \end{equation}
  Now the claim follows using the inequalities above. 
\end{proof}

\begin{lemma}\label{lem:lagr-control}
  There exists a universal constant $0<\varepsilon_1\leq\varepsilon_0$ where $\varepsilon_0$ is as in \Cref{prop:en-est} with the following property. In the setting of \Cref{lem:lag}, let $\rho>0$ be such that
  \begin{equation}
    \scurv(t,\rho)\leq\varepsilon<\varepsilon_1\quad\text{for all $t\in[0,T)$}.
  \end{equation}
  Then, for all $t\in[0,T)$,
  \begin{align}
    \int_0^t\bigl(\lambda_1^2+|\lambda_2|^{\frac43}+c_0^4\bigr)\dd\tau \leq C(\sigma,\eta,\genus\Sigma)(1+c_0^2A_0)^{\frac{16}{3}} \Big[\H_{c_0}(f(0))-\H_{c_0}(f(t)) + \frac{\sqrt{t}}{\rho^2}\Big(1+\frac{\sqrt{t}}{\rho^2}\Big)\Big]\label{eq:lagr-control}.
  \end{align}
\end{lemma}
\begin{proof}
 First, summing up the local bounds of \Cref{cor:en-est} as in \cite[Proposition 4.2]{rupp2023}, one obtains
  \begin{align}
    \int_0^t\int_{\Sigma}|A|^6\dd\mu\dd\tau &\leq C\int_{\Sigma}|A|^2\dd\mu\Big|_{t=0} + \frac{C}{\rho^4}\int_0^t\int_{\Sigma}|A|^2\dd\mu\dd\tau \\
    &\quad+ C\int_0^t\bigl(\lambda_1^2+|\lambda_2|^{\frac43}+c_0^4\bigr)\int_{\Sigma}|A|^2\dd\mu\dd\tau.
  \end{align}
  Using 
  \eqref{eq:4pi-durch-sigma}, \eqref{eq:will-vs-helf}, and the energy decay \Cref{rem:energy-decay}, we estimate
  \begin{equation}\label{eq:int|A|^2}
    \int_{\Sigma}|A|^2\dd\mu=4\W(f)-8\pi(1-\genus\Sigma)\le C(\genus\Sigma)(1+c_0^2A_0) 
  \end{equation}
  so that
  \begin{align}
    \int_0^t\int_{\Sigma}|A|^6\dd\mu\dd\tau &\leq C(\genus\Sigma)(1+c_0^2A_0)\Bigl( 1+\frac{t}{\rho^4} + \int_0^t\bigl(\lambda_1^2+|\lambda_2|^{\frac43}+c_0^4\bigr)\dd\tau\Bigr).\label{eq:lag-prop-1}
  \end{align}
  Using \Cref{lem:lag} and \eqref{eq:int|A|^2}, we infer for $C=C(\sigma,\eta,\genus\Sigma)$
  \begin{align} 
    \int_0^t \lambda_1^2\dd\tau &\leq  C(1+c_0^2A_0)\Big(\int_0^t\int_{\Sigma}|\partial_tf|^2\dd\mu\dd\tau + t\frac{c_0^2}{A_0} + tc_0^4 + \frac{1}{A_0}\int_0^t\int_{\Sigma}|A|^4\dd\mu \dd\tau \Big) \\ \label{eq:lag-prop-2}
    &\leq C (1+c_0^2A_0)\Big(\int_0^t\int_{\Sigma}|\partial_tf|^2\dd\mu\dd\tau + t\frac{c_0^2}{A_0}(1+c_0^2A_0)  \\
    &\qquad \quad + \frac{\sqrt{t(1+c_0^2A_0)}}{A_0}\Big[\int_0^t\int_{\Sigma}|A|^6\dd\mu \dd\tau\Big]^{\frac12}\Big) 
  \end{align}
  Combined with \eqref{eq:lag-prop-1}, 
  one finds
  \begin{align}
     \int_0^t\lambda_1^2\dd\tau \leq C (1+c_0^2A_0)^2&\Bigl( \int_0^t\int_{\Sigma}|\partial_tf|^2\dd\mu\dd\tau + t\frac{c_0^2}{A_0}\\
     &\quad+ \frac{\sqrt{t}}{A_0} \Bigl[1+\frac{\sqrt{t}}{\rho^2}+\Bigl(\int_0^t\bigl(\lambda_1^2+|\lambda_2|^{\frac43}+c_0^4\bigr)\dd\tau\Bigr)^{\frac12}\Bigr]\Bigr).
  \end{align}
  Using again \Cref{lem:lag}, one proceeds similarly for $\int_0^t|\lambda_2|^{\frac43}\dd\tau$ using the following observations. First,
  \begin{align} 
     \frac{1}{A_0^{\frac43}} \int_0^t\Bigl(\int_{\Sigma}|\partial_tf|\dd\mu\Bigr)^{\frac43}\dd\tau &\leq \frac{1}{A_0^{\frac43}} \int_0^t \Bigl(\int_{\Sigma}|\partial_tf|^2\dd\mu\Bigr)^{\frac23} A_0^{\frac{2}{3}}\dd\tau \leq \int_0^t\int_{\Sigma}|\partial_tf|^2\dd\mu\dd\tau + C\frac{t}{A_0^2}.
  \end{align}
  Furthermore, 
  \begin{equation}
     \frac{1}{A_0^{\frac43}}\int_0^t \Bigl(\int_{\Sigma}|A|^3\dd\mu\Bigr)^{\frac43}\dd\tau \leq \frac{1}{A_0}\int_0^t\int_{\Sigma}|A|^4\dd\mu\dd\tau,
  \end{equation}
  a term which had already been estimated in \eqref{eq:lag-prop-2}. Similarly,
  \begin{align}
     \Big(\frac{|c_0|(1+\sqrt{c_0^2A_0})}{A_0} + \frac{c_0^2}{A_0}\int_\Sigma |A|\dd\mu\Big)^\frac{4}{3} &\le C \Big(\frac{|c_0|(1+\sqrt{c_0^2A_0})}{A_0}\Big)^\frac{4}{3} + C\frac{|c_0|^\frac{8}{3}}{A_0^\frac{2}{3}}\Big(\int_\Sigma|A|^2\dd\mu\Big)^\frac{2}{3} \\
     &\le C\frac{(1+c_0^2A_0)^2}{A_0^2}.
  \end{align}%
  We thus arrive at 
  \begin{align}
      \int_0^t\bigl(\lambda_1^2 +|\lambda_2|^\frac{4}{3} +c_0^4\bigr)\dd\tau & \le C\cdot(1+c_0^2A_0)^{\frac{8}{3}}\Big(\int_0^t\int_\Sigma|\partial_tf|^2\dd\mu\dd\tau +\frac{t}{A_0^2} \\
      & \qquad + \frac{\sqrt{t}}{A_0} \Bigl[1+\frac{\sqrt{t}}{\rho^2}+\Bigl(\int_0^t\bigl(\lambda_1^2+|\lambda_2|^{\frac43}+c_0^4\bigr)\dd\tau\Bigr)^{\frac12}\Bigr]\Big).
  \end{align}
  Finally, by Simon's monotonicity formula, one has the estimate
  \begin{equation}
     \rho^2\leq C \A(f) = C A_0
  \end{equation}
  if one takes $\varepsilon_1>0$ sufficiently small, cf. \cite[Lemma 4.1]{rupp2023}. Altogether, the claim follows from \Cref{rem:energy-decay} after absorbing.
\end{proof}

\section{Construction of a blow-up}\label{sec:blowup}

\subsection{Bounds on higher-order derivatives of the second fundamental form}

For the life-span theorem and the construction of a blow-up sequence and a concentration limit, 
we also need an analogue of \Cref{cor:en-est} for higher-order derivatives of the second fundamental form $A$. This is obtained with similar arguments via local energy estimates, cf. \cite[Appendix B]{rupp2024}. This section's main result \Cref{prop:ho-en-est} is an analog of \cite[Proposition 3.7]{rupp2024} which in turn is based on \cite[Theorem 3.5]{kuwertschaetzle2001}.

\begin{proposition}\label{prop:ho-en-est:pre}
    Let $A_0,V_0>0$ satisfy \eqref{eq:iso}, $c_0\in\R$, $f\colon [0,T)\times\Sigma\to\R^3$ be an $(A_0,V_0,c_0)$-Helfrich flow, and consider $\gamma$ as in \eqref{eq:gamma}. For $m\in\N_0$, $s\ge 2m+4$, and $\phi=\nabla^m A$, one has
    \begin{align}
        &\frac{\mathrm d}{\mathrm dt}\int_\Sigma|\phi|^2\gamma^s\dd\mu + \frac12\int_\Sigma|\nabla^2 \phi|^2\gamma^s\dd\mu \le C\big(\lambda_1^2+|\lambda_2|^{\frac43}+c_0^4+\|A\|^4_{L^\infty(\{\gamma>0\})}\big)\int_\Sigma|\phi|^2\gamma^s\dd\mu \\ &\qquad \qquad + C\big(1+\lambda_1^2+|\lambda_2|^{\frac43}+c_0^4+\|A\|^4_{L^\infty(\{\gamma>0\})}\big)\int_{\{\gamma>0\}}|A|^2\dd\mu
    \end{align}
    for some $C=C(s,m,\Lambda)$ with $\Lambda$ as in \eqref{eq:gamma}.
\end{proposition}

\begin{proof}
    We will use the operator $P_r^m$ as introduced in \cite{kuwertschaetzle2002}, cf. \cite[Appendix B]{rupp2024} and \cite[Appendix A]{schlierf2024helfrich}. Writing $\partial_tf=\xi\nu$, the flow equation then reads
    \begin{equation}
        \xi = -\Delta H + P_3^0(A)+c_0P_2^0(A) + (\lambda_1+\frac12c_0^2)P_1^0(A) + \lambda_2.
    \end{equation}
    For $m\in\N_0$, a simple proof by induction using Simon's identity $\nabla^2\Delta H = \Delta^2A + P_3^2(A)$ and \cite[Lemma 2.3]{kuwertschaetzle2002} yields
    \begin{align}
        &\partial_t\nabla^mA + \Delta^2\nabla^m A = P_3^{m+2}(A)+P_5^m(A)\\ \label{eq:del_t-nablaA}
        &\quad +c_0\big(P_2^{m+2} + P_4^m(A)\big) +(\lambda_1+\frac12c_0^2)\big(P_1^{m+2}(A) + P_3^m(A)\big) + \lambda_2P_2^m(A), 
    \end{align}
    cf.\ the proof \cite[Lemma B.1]{rupp2024}. Up to the last term $\lambda_2P_2^m(A)$, we formally recover \cite[Equation (A.3)]{schlierf2024helfrich}, whereas the exceptional term $\lambda_2P_2^m(A)$ is contained in the analogous equation for the isoperimetric constrained flow \cite[Lemma B.2]{rupp2024}. Thus, the conclusion follows from \cite[Lemma 3.2]{kuwertschaetzle2002} by combining the estimates in the proofs \cite[Proposition A.1]{schlierf2024helfrich} and \cite[Propostition B.3]{rupp2024}.
\end{proof}

\begin{proposition}\label{prop:ho-en-est}
  Let $A_0,V_0>0$ satisfy \eqref{eq:iso}, $c_0\in\R$ and $\varepsilon_0$ as in \Cref{prop:en-est}. Further consider an $(A_0,V_0,c_0)$-Helfrich flow $f\colon[0,T)\times\Sigma\to\R^3$ and suppose that $\rho>0$ satisfies $T\leq T^*\rho^4$ for some $0<T^*<\infty$ and
  \begin{equation}\label{eq:kappa_le_eps}
    \scurv(t,\rho)\leq\varepsilon<\varepsilon_0\quad\text{for all $0\leq t<T$}.
  \end{equation}
  Further, assume 
  \begin{equation}\label{eq:L-mult_time_int}
    \int_0^T\bigl(\lambda_1^2+|\lambda_2|^{\frac43}+c_0^4\bigr)\dd\tau\leq\overline{L}<\infty.
  \end{equation}
  For all $t\in(0,T)$ and $m\in\N_0$, one has the local estimates
  \begin{align}
    \|\nabla^m A\|_{L^2(B_{\rho/8}(x))}&\leq C(m,T^*,\overline{L})\sqrt{\varepsilon}t^{-\frac{m}{4}},\\
    \|\nabla^m A\|_{L^{\infty}}&\leq C(m,T^*,\overline{L})\sqrt{\varepsilon}t^{-\frac{m+1}{4}},\label{eq:ho-en-est-loc}
  \end{align}
  and the full $L^2(\dd\mu)$-bounds
  \begin{equation}
    \|\nabla^m A\|_{L^2(\dd\mu)}\leq C(m,T^*,\overline{L})t^{-\frac{m}{4}}\Bigl(\int_{\Sigma}|A|^2\dd\mu\big|_{t=0}\Bigr)^{\frac12}.\label{eq:ho-en-est-full}
  \end{equation}
\end{proposition}

\begin{proof}
    Firstly, notice that by \Cref{lem:par-scal}, the time integral of
    \begin{equation}
        L(t)\vcentcolon=\lambda_1^2 + |\lambda_2|^{\frac43} + c_0^4
    \end{equation}
    is invariant under parabolic rescaling. Whence it is enough to consider the case $\rho=1$. Fix $x\in \R^3$ and define $K(t)\vcentcolon=\int_{B_1(x)}|A|^2\dd\mu$. For $0<t<T$, Equation \eqref{eq:kappa_le_eps} and \Cref{cor:en-est} imply
    \begin{equation}
        \int_0^t\int_{B_{1/2}(x)}|\nabla^2 A|^2\dd\mu\dd\tau\le C\Big(K(0) +\int_0^t(1+L)K\dd\tau\Big).
    \end{equation}
    By combining the interpolation estimates in \cite[Lemma 2.8]{kuwertschaetzle2001} and \cite[Lemma 4.2]{kuwertschaetzle2002}, and using 
    \eqref{eq:L-mult_time_int} and $T\le T^*$ it follows
    \begin{equation}
        \int_0^t\|A\|_{L^\infty(B_{1/4}(x))}^4\dd\tau \le C\Big(K(0) +\int_0^t(1+L)K\dd\tau\Big)\le C(T^*,\bar L).
    \end{equation}
    Summarizing, we have
    \begin{equation}
        \int_0^t\big(1+L+\|A\|_{L^\infty(B_{1/4}(x))}^4\big)\dd\tau \le C(T^*,\bar L)
    \end{equation}
    which is \cite[(3.18)]{rupp2024}. Now, using \Cref{prop:ho-en-est:pre} instead of \cite[Proposition B.3]{rupp2024}, the same proof of \cite[Proposition 3.7]{rupp2024} applies.
\end{proof}

\subsection{Life-span result}


The fundamental underlying result for the blow-up construction in the subsequent sections is a suitable life-span result: If one has control of the \emph{curvature concentration} and of the nonlocal Lagrange multipliers, then one obtains a (properly scaling, universal) lower bound on the maximal existence time and some control of the curvature concentration forward in time. Such a life-span result was also the starting point of the study of the Willmore flow in \cite{kuwertschaetzle2002}. Also see \cite[Propositions~5.1~and~5.2]{rupp2024} for analogous results in the context of a nonlocal flow.

\begin{proposition}\label{prop:life}
  Let $\varepsilon_1>0$ be as in \Cref{lem:lagr-control}. Then there exist universal constants $0<\overline{\varepsilon}<\varepsilon_1$ and $\overline{\delta}>0$ such that the following is satisfied. Let $A_0,V_0>0$ satisfy \eqref{eq:iso}, $c_0\in\R$ and suppose that $f_0\colon\Sigma\to\R^3$ is an embedding with $\V(f_0)=V_0$, $\A(f_0)=A_0$, and 
  \begin{equation}\label{eq:life-prop-0}
    \H_{c_0}(f_0) < \Omega(A_0,V_0,c_0),
  \end{equation}
  cf.\ \Cref{def:Omega}.
  Let $f\colon [0,T)\times\Sigma\to\R^3$ be a maximal $(A_0,V_0,c_0)$-Helfrich flow with initial datum $f_0$. If
  \begin{enumerate}[(a)]
    \item $\scurv(0,\rho)\leq\varepsilon<\overline{\varepsilon}$ for some $\rho>0$,
    \item there exists $\overline{\omega}>0$ such that, for any $t_0\in[0,\min\{T,\rho^4\overline{\omega}\}]$ with $\scurv(t,\rho)<\varepsilon_1$ for all $0\leq t<t_0$, we have $\int_0^{t_0}\bigl(\lambda_1^2+|\lambda_2|^{\frac43}+c_0^4\bigr)\dd\tau\leq\overline{\delta}$,
  \end{enumerate}
  then the maximal existence time $T$ of $f$ satisfies $T>\hat{c}\rho^4$ for some $\hat{c}=\hat{c}(\overline{\omega})\in(0,1)$,
  \begin{equation}\label{eq:lifespan-scurv-lambda-estimate}
    \scurv(t,\rho)\leq\hat{c}^{-1}\varepsilon\quad \text{for all $t\in\bigl[0,\hat{c}\rho^4\bigr]$,}\quad\text{ and }\quad \int_0^{\hat{c}\rho^4}\bigl(\lambda_1^2+|\lambda_2|^{\frac43}+c_0^4\bigr)\dd\tau\leq\overline{\delta}.
  \end{equation}
\end{proposition}
\begin{proof}
  After scaling as in \Cref{lem:par-scal}, we may assume without loss of generality that $\rho=1$. Denoting by $\Gamma>1$ the minimal number of balls of radius $\frac12$ necessary to cover $B_1(0)\subseteq\R^3$, one clearly has
  \begin{equation}\label{eq:life-prop-1}
    \scurv(t,1)\leq \Gamma \cdot\scurv(t,\frac{1}{2}).
  \end{equation}
  Then choose $\overline{\varepsilon}=\frac{\varepsilon_1}{3\Gamma}$. The constant $\overline{\delta}$ is to be chosen later.

  Note that $\scurv(t)\vcentcolon=\scurv(t,1)$ is a continuous function with $\scurv(0)\leq\varepsilon<\overline{\varepsilon}$. Therefore, for some $\omega\in(0,\overline{\omega}]$ to be chosen later,
  \begin{equation}
    t_0\vcentcolon= \sup\{t\in[0,\min\{T,\omega\}]:\scurv\leq3\Gamma\varepsilon\text{ on $[0,t)$}\} \in (0,\min\{T,\omega\}].
  \end{equation}
  By the choice of $\overline{\varepsilon}$, we have $\scurv(t)\leq 3\Gamma\varepsilon < \varepsilon_1$ for $t\in[0,t_0)$. Thus, \Cref{cor:en-est} and (b) yield
  \begin{align}
    \int_{B_{\frac12}(x)}|A|^2\dd\mu\leq \int_{B_1(x)}|A|^2\dd\mu\Big|_{t=0} + 3 \widetilde{C}\Gamma\varepsilon t + 3\widetilde{C}\Gamma\varepsilon\overline{\delta},
  \end{align} 
  for all $0\leq t<t_0$ where $\widetilde{C}$ is the universal constant ``$\,C\,$'' from \Cref{cor:en-est}. Choose $\omega=\min\{(6\widetilde{C}\Gamma)^{-1},\overline{\omega}\}$ and $\overline{\delta}=(6\widetilde{C}\Gamma)^{-1}$. Then
  \begin{equation}
    \int_{B_{\frac12}(x)}|A|^2\dd\mu\leq \int_{B_1(x)}|A|^2\dd\mu\Big|_{t=0} + \frac{\varepsilon}{2}\frac{1}{\omega}t+\frac{\varepsilon}{2} \leq 2\varepsilon\quad\text{for all $0\leq t<t_0$}.
  \end{equation}
  Using \eqref{eq:life-prop-1}, if $t_0<\min\{T,\omega\}$, one finds $\scurv(t)\leq 2\Gamma\varepsilon$ for all $0\leq t<t_0$. However, by the definition of $t_0$ and continuity of $\scurv$, in this case, one also has $\scurv(t_0)=3\Gamma\varepsilon$, a contradiction!

  \textbf{Assumption: } $t_0=T\le\omega$. Using the specific choices for $t_0$ and $\overline{\varepsilon}$ made above, one again obtains $\scurv(t)\leq 3\Gamma\varepsilon<\varepsilon_1$ for $0\leq t<T$. Particularly, using $T\leq\overline{\omega}$, Assumption (b) yields $\int_0^T\bigl(\lambda_1^2+|\lambda_2|^{\frac43}+c_0^4\bigr)\dd\tau\leq\overline{\delta}$. So \Cref{prop:ho-en-est} implies for $0<\xi<T$ and $m\in\N_0$
  \begin{align}
    \|\nabla^m A\|_{L^{\infty}}  &\leq C(m,\overline{\omega},\overline{\delta},\xi)\quad\text{and}\\
    \|\nabla^m A\|_{L^{2}(\dd\mu)}  &\leq C(m,\overline{\omega},\overline{\delta},\xi,\operatorname{genus}\Sigma,\H_{c_0}(f_0),c_0^2A_0)
  \end{align}
  for $t\in[\xi,T)$, respectively. Thus, for $t\in [\xi,T)$, using \Cref{prop:cmc-control} combined with \eqref{eq:life-prop-0} 
  to control the denominators in \eqref{eq:def-l1}, \eqref{eq:def-l2} as in \eqref{eq:denom-cmc-estimate},
  \begin{align}
    |\lambda_1|^2&\leq C(\genus\Sigma, A_0,V_0,c_0,\Omega(A_0,V_0,c_0)-\H_{c_0}(f_0))\ \|\nabla\H_{c_0}(f)\|_{L^2(\dd\mu)}^2 \\
    &\leq C(\genus\Sigma, A_0,V_0,c_0,\Omega(A_0,V_0,c_0)-\H_{c_0}(f_0),\xi,\overline{\omega},\overline{\delta})
  \end{align}
  and similarly also for $\lambda_2$.  With the same arguments as in \cite[pp. 330 -- 332]{kuwertschaetzle2002}, one deduces that $f(t)$ converges smoothly to an immersion $f(T)$ as $t\nearrow T$. Thus, one can restart the flow at time $T$ which contradicts the maximality of $T$. Therefore, $T>\omega=t_0$.  Consequently, as observed above for $t< t_0 =\omega<T$, we have $\kappa(t)\leq 2\Gamma\varepsilon <\frac{2}{3}\varepsilon_1$. Choosing $\hat{c}\vcentcolon=\min\{\omega,(2\Gamma)^{-1}\}>0$, the estimates in \eqref{eq:lifespan-scurv-lambda-estimate} hence follow.
\end{proof}

Observe that, for the life-span result \Cref{prop:life}, we assume some sort of boundedness of the nonlocal Lagrange multipliers. More precisely, the energy threshold \eqref{eq:life-prop-0} is needed to bound their denominator, see \eqref{eq:denom-cmc-estimate} and \Cref{prop:cmc-control}, and the $L^p(\dd t)$ control assumed in \Cref{prop:life}(b) is needed in the localized estimates of \Cref{cor:en-est}. The following corollary now applies \Cref{lem:lagr-control} to obtain \Cref{prop:life}(b) under the second part of the energy assumption in \eqref{eq:life-cor-0}.

\begin{corollary}\label{cor:life}
    Let $\overline{\varepsilon}>0$ be as in \Cref{prop:life},  $A_0,V_0>0$ satisfy \eqref{eq:iso}, $c_0\in\R$, and $f_0\colon\Sigma\to\R^3$ be an embedding with $\A(f_0)=A_0$, $\V(f_0)=V_0$, satisfying 
    \begin{equation}\label{eq:life-cor-0}
    \sqrt{\H_{c_0}(f_0)} \leq \min\Bigl\{\sqrt{\Omega(A_0,V_0,c_0)(1-\eta)}, \Bigl(\sqrt{\frac{4\pi}{\sigma}}(1-\eta) - \frac{c_0}2\sqrt{A_0}\Big)\Bigr\}
  \end{equation}
    for some $\eta\in(0,1)$. If $f\colon[0,T)\times\Sigma\to\R^3$ is a maximal $(A_0,V_0,c_0)$-Helfrich flow starting in $f_0$ with
    \begin{enumerate}[(i)]
        \item $\scurv(0,\rho)\leq\varepsilon<\overline{\varepsilon}$ for some $\rho>0$ and
        \item $\H_{c_0}(f_0)-\lim_{t\nearrow T}\H_{c_0}(f(t))\leq \overline{d}$ where $\overline{d}=\overline{d}(\eta,\sigma,\genus\Sigma,c_0^2A_0)>0$,
    \end{enumerate}
    then $T>\hat{c}\rho^4$ where $\hat{c}=\hat{c}(\eta,\sigma,\genus\Sigma,c_0^2A_0)\in(0,1)$ and we have the estimates \eqref{eq:lifespan-scurv-lambda-estimate}.
\end{corollary}
\begin{proof}
    One only needs to check that the assumptions (a) and (b) in \Cref{prop:life} are satisfied, for some appropriately chosen $\overline{d}$ which will be determined by \Cref{lem:lagr-control}. Clearly, (a) is satisfied by (i). For (b), $\overline{\omega}>0$ is to be chosen later and suppose that $0<t_0\leq \min\{T,\rho^4\overline{\omega}\}$ satisfies $\scurv(t,\rho)<\varepsilon_1$ for all $0\leq t < t_0$. By the energy threshold \eqref{eq:life-cor-0}, \Cref{lem:lagr-control} applies and \eqref{eq:lagr-control} together with (ii) yields
    \begin{equation}\label{eq:cor-life-1}
        \int_0^{t_0} \big(\lambda_1^2+|\lambda_2|^{\frac43}+c_0^4\big)\dd t \leq C(\eta,\sigma,\genus\Sigma,c_0^2A_0)\big(\overline{d} + \overline{\omega}^{\frac12}+\overline{\omega}\big) \leq\overline{\delta}
    \end{equation}    
    if one takes $\overline{d}=\overline{d}(\eta,\sigma,\genus\Sigma,c_0^2A_0)>0$ and $\overline{\omega}=\overline{\omega}(\eta,\sigma,\genus\Sigma,c_0^2A_0)>0$ sufficiently small. 
\end{proof}

\subsection{Existence of a blow-up and resulting concentration limit}

\begin{lemma}\label{lem:ex-blow-up}
     Let $\overline{\varepsilon}>0$ be as in \Cref{prop:life},  $A_0,V_0>0$ satisfy \eqref{eq:iso}, $c_0\in\R$, and let $f_0\colon\Sigma\to\R^3$ be an embedding with $\V(f_0)=V_0$, $\A(f_0)=A_0$, satisfying \eqref{eq:life-cor-0} for some $\eta\in(0,1)$. Fix $\hat{c}=\hat{c}(\eta,\sigma,\genus\Sigma,c_0^2A_0)$ as in \Cref{cor:life} and consider a maximal $(A_0,V_0,c_0)$-Helfrich flow $f\colon[0,T)\times\Sigma\to\R^3$ starting in $f_0$.

     Then there exist times $t_j\nearrow T$, radii $(r_j)_{j\in\N}\subseteq(0,\infty)$ and points $(x_j)_{j\in\N}\subseteq\R^3$ such that each $f_j\colon[0,(T-t_j)/r_j^4)\times\Sigma\to\R^3$,
     \begin{equation}
         f_j(t,p)\vcentcolon=\frac{1}{r_j}\big(f(t_j+r_j^4t,p)-x_j\big)
     \end{equation}
     is a maximal $((1/r_j^2)A_0,(1/r_j^3)V_0,r_jc_0)$-Helfrich flow, $r_j\leq C\big(\sqrt{A_0\H_{c_0}(f_0)} + \frac12 (c_0)_+A_0\big)$ for a universal constant $C>0$ where $(c_0)_+=\max\{c_0,0\}$, and
     \begin{enumerate}[(i)]
         \item \label{it:blow_up1} $t_j+r_j^4\hat{c}<T$;
         \item \label{it:blow_up2} $\scurv_j(t,1)\leq\overline{\varepsilon}$ for all $0\leq t\leq \hat{c}$ and $\int_0^{\hat{c}}\bigl(\lambda_1(f_j)^2+|\lambda_2(f_j)|^{\frac43}+(r_jc_0)^4\bigr)\dd\tau\leq\overline{\delta}$;
         \item \label{it:blow_up3}  $\inf_{j\in\N}\int_{B_1(0)}|A_{f_j(\hat{c})}|^2\dd\mu_{f_j(\hat{c})}\ge \overline\alpha$ for some $\overline \alpha>0$ depending only on $\overline\varepsilon$ and $\hat c$.
     \end{enumerate}
\end{lemma}
\begin{proof}
  The same proof as \cite[Lemma 6.6 in Article C]{rupp2022PhD} gives the existence of $\alpha>0$ and radii $(r_t)_{0\leq t<T}\subseteq(0,\infty)$ with
  \begin{equation}\label{eq:lem-ex-bu-1}
    \alpha <\scurv(t,r_t)<\hat{c}\overline{\varepsilon}.
  \end{equation}
  Indeed, these $r_t$ satisfy $r_t \le \diam f(t,\Sigma)$ and thus, arguing as in \eqref{eq:diam-bd-c0geq0} and \eqref{eq:diam-bd-c0leq0}, 
  \begin{align}
    r_t&\leq 
    C \begin{cases}
        \sqrt{A_0\W(f(t))} &\text{if $c_0\geq 0$}\\
        \sqrt{A_0\H_{c_0}(f(t))}&\text{if $c_0<0$}
    \end{cases} \leq C\sqrt{A_0}\begin{cases}
        \sqrt{\H_{c_0}(f_0)} + \frac12\sqrt{c_0^2 A_0}&\text{if $c_0\geq 0$}\\
        \sqrt{\H_{c_0}(f_0)}&\text{if $c_0<0$}
    \end{cases}\\
    &= C\big(\sqrt{A_0\H_{c_0}(f_0)} + \frac12 (c_0)_+A_0\big),
  \end{align}
  using \Cref{rem:energy-decay} and \eqref{eq:will-vs-helf} in the case where $c_0\geq 0$. Moreover, there exists $t_0\in[0,T)$ such that, for $t_0\leq t<T$, 
  \begin{equation}\label{eq:lem-ex-bu-2}
    \H_{c_0}(f(t))-\lim_{\tau\nearrow T}\H_{c_0}(f(\tau))\leq\overline{d}
  \end{equation}
  with $\overline{d}$ as in \Cref{cor:life}. Since both constants $\overline d$ and $\hat c$ are invariant under parabolic rescaling and using the uniform energy bound in \eqref{eq:life-cor-0}, \Cref{cor:life} implies $t+r_t^4\hat c < T$. In particular, Claim \eqref{it:blow_up1} will be satisfied. Proceeding as in the proof of \cite[Lemma 5.5]{schlierf2024helfrich}, one deduces from the upper bound on $r_t$ that there exist $t_0\leq t_j\nearrow T$ and $0<M<2$ with 
  \begin{equation}
    r_{t_j+r_j^4\hat{c}} \leq Mr_{t_j}\quad\text{for all $j\in\N$}.
  \end{equation}
  Write $r_j\vcentcolon=r_{t_j}$. Claim \eqref{it:blow_up2} follows from \eqref{eq:lifespan-scurv-lambda-estimate}. 
  By a covering argument, there is $N=N(M)>0$ with
  \begin{equation}
    \scurv(t,Mr)\leq N\scurv(t,r)\quad\text{for all $0\leq t<T$ and $r>0$}.
  \end{equation}
  Then $\scurv(t_j+r_j^4\hat{c},r_j)\geq\frac1N\scurv(t_j+r_j^4\hat{c},Mr_j)\geq\frac1N\scurv(t_j+r_j^4\hat{c},r_{t_j+r_j^4\hat{c}})\geq\frac1N\alpha$, using \eqref{eq:lem-ex-bu-1}. Appropriately choosing $x_j\in\R^3$ yields \eqref{it:blow_up3}.
\end{proof}

Using the a priori control of $\lambda_1,\lambda_2$ obtained in \Cref{sec:lagr}, one finds

\begin{proposition}[Existence and properties of a concentration limit]\label{prop:concentration-limit}
    In the setting of \Cref{lem:ex-blow-up}, there exists a complete, orientable surface $\hat{\Sigma}\neq \emptyset$ without boundary and a proper immersion $\hat{f}\colon\hat{\Sigma}\to\R^3$ such that, after passing to a subsequence, $r_j\to r\in[0,\infty)$ and, for a suitable choice of orientation on $\hat{\Sigma}$,
    \begin{enumerate}[(i)]
        \item \label{it:blow_up21} $\hat{f}_j\vcentcolon=f_j(\hat{c})\to\hat{f}$ smoothly on compact subsets after positive reparametrizations;
        \item \label{it:blow_up22} $\int_{\bar B_1(0)}|\hat{A}|^2\dd\hat{\mu}>0$;
        \item \label{it:blow_up23} $\hat{f}$ is a constrained $\H_{rc_0}$-Helfrich immersion;
        \item \label{it:blow_up24} $\A(\hat{f}_j)\to\infty$ if and only if $r=0$. In this case, $\hat{f}$ is a Willmore immersion.
    \end{enumerate}
\end{proposition}
\begin{proof} 
  As $r_j\leq C\big(\sqrt{A_0\H_{c_0}(f_0)} + \frac12 (c_0)_+A_0\big)$ by \Cref{lem:ex-blow-up}, after passing to a subsequence, we have $r_j\to r$ with $0\leq r<\infty$. Using \eqref{it:blow_up2} in \Cref{lem:ex-blow-up} and \Cref{prop:ho-en-est}, one finds for $0 < t\leq \hat{c}$ and, for all $m\in\N_0$,
  \begin{align} \label{eq:prop-conc-lim-1}
    \|\nabla^m A_j\|_{L^{\infty}}&\leq C(m,\hat{c},\overline{\delta})t^{-\frac14(m+1)},\\ \label{eq:prop-conc-lim-2}
    \|\nabla^m A_j\|_{L^{2}(\dd\mu_j)}&\leq C(m,\hat{c},\overline{\delta},\genus\Sigma,\H_{c_0}(f_0),c_0^2A_0)t^{-\frac{m}{4}}
  \end{align}
  using 
  \begin{align}
    \int_{\Sigma}|A_j|^2\dd\mu_j\Big|_{t=0}& = 4\mathcal W(f(t_j)) -4\pi\chi(\Sigma) \\
    &\le C(\genus\Sigma)(1 + \mathcal W(f(t_j))) \le C(\genus\Sigma)(1+\H_{c_0}(f_0) + c_0^2A_0).        
  \end{align}
  Similarly, by Simon's monotonicity formula 
  \cite[Equation (1.3)]{simon1993},
  \begin{equation}
      \mu_{\hat f_j}(B_R(0))\le CR^2\mathcal W(f(t_j+\hat c r_j^4))\le R^2C(\H_{c_0}(f_0),c_0^2A_0).
  \end{equation}
  The localized version of Langer's compactness theorem in \cite[Theorem 4.2]{kuwertschaetzle2001}, also cf. \cite[Appendix A]{rupp2023}, yields the following. There exists a complete $2$-manifold $\hat{\Sigma}$ without boundary and a proper immersion $\hat{f}\colon\hat{\Sigma}\to\R^3$ such that, after passing to a subsequence, $\hat{f}_j=f_j(\hat{c})\to \hat{f}$ smoothly on compact subsets of $\R^3$ after reparametrization. That is, writing $\hat{\Sigma}(R)=\{p\in\hat\Sigma\mid|\hat{f}(p)|<R\}$ for $R>0$, one has a representation
  \begin{equation}\label{eq:conv-on-cpct-subsets}
    \hat{f}_j\circ\phi_j = \hat{f} + u_j\quad\text{on $\hat{\Sigma}(j)$}
  \end{equation}
  where
  \begin{align}
    &\phi_j\colon\hat{\Sigma}(j)\to U_j\subseteq\Sigma\text{ is a diffeomorphism},\\
    &\{|f_j|<R\} \subseteq \phi_j(\hat{\Sigma}(j))\text{ for $j\geq j(R)$, for all $R>0$},\\
    &u_j\in C^{\infty}(\hat{\Sigma}(j),\R^3) \text{ is normal along $\hat{f}$},\\
    &\|\hat{\nabla}^mu_j\|_{L^{\infty}(\hat{\Sigma}(j))}\to 0\text{ for all $m\in\N_0$}.
  \end{align} 
  This already proves \eqref{it:blow_up21}. Moreover, \eqref{it:blow_up22} simply follows from \Cref{lem:ex-blow-up}\eqref{it:blow_up3} which also gives $\hat{\Sigma}\neq\emptyset$. 
  
  Now fix some $0<\xi<\hat{c}$ and consider the flows
  \begin{equation}
    \tilde{f}_j\colon[\xi,\hat{c}]\times\hat{\Sigma}(j)\to\R^3,\quad\tilde{f}_j(t,p)=\frac{1}{r_j}(f(t_j+r_j^4t,\phi_j(p))-x_j).
  \end{equation}
  Particularly, $\tilde{f}_j(\hat{c})=\hat{f}+u_j$ smoothly converges to $\hat{f}$ on compact subsets of $\hat{\Sigma}$. Also, the curvature estimates in \eqref{eq:prop-conc-lim-1} remain valid for the reparametrized flows $\tilde{f}_j$ where the powers of $t^{-1}$ are uniformly controlled by powers of $\xi^{-1}$. In view of \Cref{prop:scaling_Omega}, we may combine Equations \eqref{eq:def-l1}--\eqref{eq:denom-cmc-estimate} with \Cref{prop:cmc-control} to deduce that for some constant $C$ depending on $A_0,V_0,c_0,\eta,\genus\Sigma$, 
  there holds
  \begin{align}
    |\lambda_1(f_j)|&\leq C \|\nabla\H_{r_jc_0}({f}_j)\|_{L^2(\Sigma,\dd{\mu}_j)},\\
    |\lambda_2({f}_j)|&\leq \frac{C}{\sqrt{\A({f}_j)}} \|\nabla\H_{r_jc_0}({f}_j)\|_{L^2(\Sigma,\dd{\mu}_j)} = r_j \frac{C}{\sqrt{A_0}} \|\nabla\H_{r_jc_0}({f}_j)\|_{L^2(\Sigma,\dd{\mu}_j)}.\label{eq:conc-limit-bound-on-lambda2}
  \end{align}
  Using \eqref{eq:prop-conc-lim-1} as well as $r_j\leq C\big(\sqrt{A_0\H_{c_0}(f_0)} + \frac12 (c_0)_+A_0\big)$, it follows
  \begin{equation} \label{eq:blow-up:control_lambdas}
      |\lambda_1(f_j)|+|\lambda_2(f_j)|\le C(\hat{c},\overline{\delta},\genus\Sigma,\xi,c_0,A_0,V_0,\eta).
  \end{equation}
  Notice that by the same proof, \eqref{eq:del_t-nablaA} remains valid if on the left hand side, the second fundamental form $A$ is replaced by the mean curvature $H$. Combining this with \eqref{eq:prop-conc-lim-1}, \eqref{eq:prop-conc-lim-2}, \eqref{eq:blow-up:control_lambdas}, one finds for $\xi\leq t\leq \hat{c}$ that
  \begin{align}
    \|\partial_t \nabla^m A_j\|_{L^{\infty}}+\|\partial_t \nabla^m A_j\|_{L^{2}(\dd\mu_j)} &\leq  C,\\
    \|\partial_t \nabla^m H_j\|_{L^{\infty}} +\|\partial_t \nabla^m H_j\|_{L^{2}(\dd\mu_j)} &\leq  C
  \end{align}
  for a constant $C=C(m,\hat{c},\overline{\delta},\genus\Sigma,\xi,c_0,A_0,V_0,\eta)$. 
  Therefore, using \eqref{eq:def-l1}, \eqref{eq:def-l2}, \eqref{eq:denom-cmc-estimate}, \Cref{prop:cmc-control}, and \eqref{eq:prop-conc-lim-1}, \eqref{eq:prop-conc-lim-2}, we infer for $\xi\leq t\leq\hat{c}$
  \begin{align}
    |\partial_t \lambda_1({f}_j)| +|\partial_t \lambda_2({f}_j)| \leq  C.
  \end{align}
  for a constant $C=C(m,\hat{c},\overline{\delta},\genus\Sigma,\xi,c_0,A_0,V_0,\eta)$.
  Therefore, arguing as in \cite[Article C, pp. 25--26]{rupp2022PhD}, after passing to a subsequence, the flows $(\tilde{f}_j)_{j\in\N}$ converge in $C^1([\xi,\hat{c}],C^{\infty}(P))$ for any compact $P\subseteq\hat{\Sigma}$ to a limiting flow $\tilde{f}_{\mathrm{lim}}\colon[\xi,\hat{c}]\times\hat{\Sigma}\to\R^3$. With similar arguments, for the unit-normal fields, one may assume $\nu_j(t,\phi_j(x)) \to\tilde{\nu}_{\mathrm{lim}}(t,x)$ in $C^1([\xi,\hat{c}],C^{\infty}(P))$ where $\tilde{\nu}_{\mathrm{lim}}(t,\cdot)$ is a smooth unit normal field along $\tilde{f}_{\mathrm{lim}}(t)$. Particularly, $\hat{\Sigma}$ is again orientable and we can fix an orientation on $\hat{\Sigma}$ such that $\tilde{\nu}_{\mathrm{lim}}$ is induced by $\tilde{f}_{\mathrm{lim}}$ via \eqref{eq:def-nu}. In particular, each $\phi_j\colon\hat{\Sigma}(j)\to U_j\subseteq \Sigma$ is positive for the induced orientations on $\hat{\Sigma}(j)\subseteq\hat{\Sigma}$ and $U_j\subseteq\Sigma$, respectively. Further, we may assume
  \begin{equation}
    \lambda_1({f}_j) \to \lambda_{1,\mathrm{lim}}\quad\text{and}\quad \lambda_2({f}_j)\to\lambda_{2,\mathrm{lim}}\quad\text{in $C^0([\xi,\hat{c}],\R)$}.
  \end{equation} 
  Then, for any compact $P\subseteq\hat{\Sigma}$,
  \begin{align}
    \int_{\xi}^{\hat{c}} \int_P |\partial_t\tilde{f}_{\lim}|^2\dd\tilde{\mu}_{\lim}\dd t &= \lim_{j\to\infty} \int_{\xi}^{\hat{c}} \int_P |\partial_t\tilde{f}_j|^2\dd\tilde{\mu}_j\dd t \leq \lim_{j\to\infty} \int_{\xi}^{\hat{c}} \int_{\Sigma} |\partial_t f_j|^2\dd\mu_j\dd t\\
    &= 2\lim_{j\to\infty} \big(\H_{r_jc_0}(f_j(\xi))-\H_{r_jc_0}(f_j(\hat{c}))\big) \\
    &= 2\lim_{j\to\infty} \big(\H_{c_0}(f(t_j+r_j^4\xi))-\H_{c_0}(f(t_j+r_j^4\hat{c}))\big) = 0,
  \end{align}
  using \Cref{rem:energy-decay}, \eqref{eq:scaling_Helfrich}, and the smooth convergence established above. So $\tilde{f}_{\lim}(t)=\tilde{f}_{\lim}(\hat{c}) = \hat{f}$ on $\hat{\Sigma}$ for all $\xi\leq t\leq \hat{c}$. Setting $\hat{\lambda}_1=\lambda_{1,\lim}(\hat{c})$, $\hat{\lambda}_2=\lambda_{2,\lim}(\hat{c})$, $\hat{\nu}=\tilde{\nu}_{\lim}(\hat{c})$, the fact that, due to \Cref{lem:par-scal} and the choice of orientation on $\hat{\Sigma}$, each $\tilde{f}_j$ is an $(r_j^{-2}A_0,r_j^{-3}V_0,r_jc_0)$-Helfrich flow, and $r_j\to r$ yield with the above convergence in $C^1([\xi,\hat{c}],C^{\infty}(P))$ for any compact $P\subseteq\hat{\Sigma}$ 
  \begin{align}
    &\big(\Delta\hat{H} + |\hat{A}^0|^2\hat{H} -rc_0(|\hat{A}^0|^2-\frac12\hat{H}^2) -(\hat{\lambda}_1+\frac12r^2c_0^2)\hat{H} - \hat{\lambda}_2\big)\hat{\nu}\\
    &\quad =-\lim_{j\to\infty} \partial_t\tilde{f}_j(\hat{c}) = -\partial_t\tilde{f}_{\lim}(\hat{c}) = 0\quad\text{on $\hat{\Sigma}$}.
  \end{align}
  So $\hat{f}$ is a constrained $\H_{rc_0}$-Helfrich immersion. This proves Claim \eqref{it:blow_up23}. 
  
  For \eqref{it:blow_up24}, since
  \begin{equation}
    \A(\hat{f}_j) = \A\Big(\frac{1}{r_j}\big(f(t_j+r_j^4\hat{c})-x_j\big)\Big) = \frac{A_0}{r_j^2},
  \end{equation}
  one finds $r=0$ if and only if $\A(\hat{f}_j)\to \infty$. With \eqref{it:blow_up23}, we only need to show that, in this case, $\hat{\lambda}_1=\hat{\lambda}_2=0$. To this end, suppose that $r=0$ or, equivalently, $\A(\hat{f}_j)\to\infty$. 
  By \Cref{lem:lag} and \Cref{rem:energy-decay}, for all $j\in\N$ and $\xi\in(0,\hat{c})$, 
  \begin{align}
    &\int_{\xi}^{\hat{c}} |\lambda_1(f_j)|^2\dd t\leq C(\sigma,\eta,\genus\Sigma)(1+c_0^2A_0) \Big(\H_{c_0}(t_j+r_j^4\xi)-\H_{c_0}(t_j+r_j^4\hat{c}) \\
    &\quad + \frac{1}{\frac{1}{r_j^2}A_0}\int_{\xi}^{\hat{c}}\Big[r_j|c_0| + r_j^2c_0^2\int_{\Sigma}|A_j|\dd\mu_j +  \int_{\Sigma}|A_j|^3\dd\mu_j \Big]^2\dd t\Big) \to 0
  \end{align}
  as $j\to\infty$, using \eqref{eq:prop-conc-lim-1},\eqref{eq:prop-conc-lim-2}, and $r_j\to 0$. So $\hat{\lambda}_1=\lim_{j\to\infty}\lambda_1(f_j(\hat{c}))=0$. Furthermore, as in \eqref{eq:conc-limit-bound-on-lambda2}, one finds 
  \begin{equation}
      |\lambda_2({f}_j)|\leq \frac{C(\sigma)}{\sqrt{\A({f}_j)}} \|\nabla\H_{r_jc_0}({f}_j)\|_{L^2(\Sigma,\dd{\mu}_j)} = r_j \frac{C(\sigma)}{\sqrt{A_0}} \|\nabla\H_{r_jc_0}({f}_j)\|_{L^2(\Sigma,\dd{\mu}_j)}
  \end{equation}
 and $\|\nabla\H_{r_jc_0}({f}_j)\|_{L^2(\Sigma,\dd{\mu}_j)}$ is uniformly bounded in $j\in\N$ for $\xi\leq t\leq \hat{c}$ by \eqref{eq:prop-conc-lim-1}--\eqref{eq:prop-conc-lim-2}. Therefore, we conclude $\hat{\lambda}_2=\lim_{j\to\infty}\lambda_2(f_j(\hat{c}))=0$. 
\end{proof}

\section{Stability via a constrained \L ojasiewicz-Simon inequality}\label{sec:loja}

In this section, we prove global existence for initial data sufficiently close to a local minimizer. 
Our method is based on a {\L}ojasiewicz--Simon gradient inequality and the arguments in \cite{chillfasangovaschaetzle2009} for the Willmore flow.
A version of this inequality for the constrained Canham--Helfrich energy has already been formally derived in a model for fluid vesicle dynamics, see \cite[Theorem 1.4]{Lengeler2018}. However, this approach works in a different Sobolev regime and requires a graphical representation close to an embedded surface, thus it is not applicable in our context. Instead, as in \cite{rupp2024}, we use the sufficient conditions of \cite{Rupp2020}.

\begin{lemma}\label{lem:Loja}
    Let $f\colon \Sigma \to \R^3$ be a constrained $\mathcal{H}_{c_0}$-Helfrich immersion with $H_f\not\equiv const$. Then there exist $C,\rho>0$ and $\theta\in (0,\frac12]$ such that for all immersions $h\in W^{4,2}(\Sigma;\R^3)$ with $\Vert h-f\Vert_{W^{4,2}}\leq \rho$, $\mathcal{A}(h)=\mathcal{A}(f)$, and $\mathcal{V}(h)=\mathcal{V}(f)$ we have
    \begin{align}
        |\mathcal{H}_{c_0}(h)-\mathcal{H}_{c_0}(f)|^{1-\theta} \leq C \Vert 2\nabla \mathcal{H}_{c_0}(h) -\lambda_1(h) H_h\nu_h - \lambda_2(h) \nu_h\Vert_{L^2(\mathrm{d}\mu_h)},
    \end{align}
    where $\lambda_1(h),\lambda_2(h)$ are as in \eqref{eq:def-l1} and \eqref{eq:def-l2}.
\end{lemma}

We will not prove this lemma here, since it can be deduced with minor modifications from the case of the Willmore energy with constrained isoperimetric ratio in \cite[Section 5.3]{rupp2024}.
Indeed, the energy $\mathcal{H}_{c_0}$ is essentially a lower order perturbation of the Willmore functional, so the necessary Fredholm properties can be deduced from the second variation of $\mathcal{W}$ computed in \cite{chillfasangovaschaetzle2009}. It remains to check the conditions on the constraints in \cite[Corollary 5.2]{Rupp2020}. Since the constraints on the area and the volume are of lower order, compared to the energy $\mathcal{H}_{c_0}$, condition (v) in \cite[Corollary 5.2]{Rupp2020} follows from the representation of the $L^2$-gradients in \Cref{prop:variations}. By a direct computation, the constraints are nondegenerate, i.e., assumption (vi) of  \cite[Corollary 5.2]{Rupp2020} is satisfied, if $H_f\not\equiv const$. 

As a consequence of \Cref{lem:Loja}, we find the following asymptotic stability result.

\begin{lemma}\label{lem:LojaAsymStabil}
		Let $\hat{f}\colon \Sigma\to\R^3$ be a constrained $\mathcal{H}_{c_0}$-Helfrich immersion with $\mathcal{A}(\hat{f})=A_0$, $\mathcal{V}(\hat{f})=V_0$, and $H_{\hat{f}}\not \equiv const$. 
		 Let $\alpha\in(0,1)$, $k\in \N$ with $k\geq 4$, and $\delta>0$. Then there exists $\varepsilon=\varepsilon(\hat{f},\alpha,k,\delta)>0$ such that if $f\colon [0,T)\times \Sigma\to\R^3$ is an $(A_0,V_0,c_0)$-Helfrich flow satisfying
		\begin{enumerate}[(i)]
			\item $\Vert {f_0-\hat{f}}\Vert_{C^{k,\alpha}}<\varepsilon$;
			\item ${\mathcal{H}_{c_0}}(f(t))\geq {\mathcal{H}_{c_0}}(\hat{f})$ whenever $\Vert {f(t)\circ \Phi(t) - \hat{f}}\Vert_{C^k}\leq \delta$, for some positive diffeomorphisms $\Phi(t) \colon\Sigma\to\Sigma$;
		\end{enumerate}
		then the flow exists globally, i.e., we may take $T=\infty$. Moreover, as $t\to\infty$, it converges smoothly after reparametrization by some positive diffeomorphisms $\tilde{\Phi}(t)\colon\Sigma\to\Sigma$ to a constrained $\mathcal{H}_{c_0}$-Helfrich immersion $f_\infty$, satisfying ${\mathcal{H}_{c_0}}(f_\infty)={\mathcal{H}_{c_0}}(\hat{f})$ and $\Vert {f_{\infty}-\hat{f}}\Vert_{C^k}\leq\delta$.
\end{lemma}

This is exactly the counterpart to \cite[Lemma~4.1]{chillfasangovaschaetzle2009}, \cite[Lemma~7.9]{rupp2023}, and \cite[Lemma~5.8]{rupp2024} for our gradient flow situation. As a consequence, $r>0$ in \Cref{prop:concentration-limit} implies convergence.

\begin{lemma}
  \label{lem:r>0 gives conv}
    Suppose that $\hat{f}\colon \hat\Sigma\to\R^3$ is a concentration limit as in \Cref{prop:concentration-limit} of a maximal $(A_0,V_0,c_0)$-Helfrich flow $f\colon [0,T)\times \Sigma\to\R^3$. If $r>0$, then $\hat{\Sigma}\cong\Sigma$ is compact. Further, if ${H}_{\hat{f}}\not\equiv const$, then  the flow exists for all times, i.e., $T=\infty$, and, as $t\to\infty$, converges smoothly after positive reparametrizations to a constrained $\mathcal{H}_{c_0}$-Helfrich immersion $f_\infty$.
\end{lemma}

\begin{proof}
    By \Cref{prop:concentration-limit}\eqref{it:blow_up24}, $r>0$ implies $\sup_j \mathcal{A}(\hat{f}_j)<\infty$. Simon's diameter estimate \cite[Lemma 1.1]{simon1993} then implies that $\hat{f}(\hat{\Sigma})$ is compact, and consequently, by \cite[Lemma~4.3]{kuwertschaetzle2001}, we have that $\hat{\Sigma}\cong \Sigma$ is compact. In fact, \cite[Lemma~4.3]{kuwertschaetzle2001} yields that, in \eqref{eq:conv-on-cpct-subsets}, one may without loss of generality arrange that $\hat{\Sigma}(j)=\hat{\Sigma}$ and $U_j=\Sigma$ such that $\phi_j\colon\hat{\Sigma}\to\Sigma$ is a positive diffeomorphism for all $j\in\N$, after passing to a subsequence. Using \eqref{eq:conv-on-cpct-subsets},
    \begin{align}
        f(t_j+\hat{c}r_j^4,\phi_j\circ\phi_1^{-1}(\cdot))-x_j = r_j \hat{f}_j\circ(\phi_j\circ\phi_1^{-1}) \to r \hat{f}\circ\phi_1^{-1}\qquad \text{smoothly on }\Sigma.
    \end{align}
    Moreover, $\phi_j\circ\phi_1^{-1}\colon\Sigma\to\Sigma$ is positive for each $j\in\N$. By \Cref{prop:concentration-limit}\eqref{it:blow_up23} and \eqref{eq:scaling_Helfrich}, we have that $r\hat{f}\circ\phi_1^{-1}$ is a constrained $\H_{c_0}$-Helfrich immersion since $\phi_1^{-1}\colon\Sigma\to\hat{\Sigma}$ is positive. Furthermore, $H_{r\hat{f}\circ\phi_1^{-1}} = r^{-1}H_{\hat{f}\circ\phi_1^{-1}} = r^{-1}H_{\hat{f}}\circ\phi_1^{-1} \not\equiv const$,
    so that for $j$ large enough, the assumptions of \Cref{lem:LojaAsymStabil} are satisfied and the second statement follows.
\end{proof}

\section{Convergence for spherical immersions}\label{sec:conv_sphere}

In view of \Cref{lem:r>0 gives conv}, we want to exclude the case $r=0$ in \Cref{prop:concentration-limit}, i.e., the case of a noncompact concentration limit. To this end, we proceed by contradiction and first deduce some properties of a noncompact concentration limit by adapting arguments in \cite[Appendix~A]{kuwertschaetzle2004} to the monotonicity formula for the Helfrich energy $\H_{c_0}$ found in \cite{ruppscharrer2023}. In the following, we define the function which is the $\H_{c_0}$-analog to \cite[Equation~(A.4)]{kuwertschaetzle2004}.

\begin{definition}
  For an immersion $h\colon\Sigma\to\R^3$ of an oriented, complete surface $\Sigma$ without boundary and $c_0\in\R$, $\rho>0$, define
  \begin{align}
      \gamma_{h,c_0}(\rho) &= \frac{\mu(\bar{B}_{\rho}(0))}{\rho^2} + \frac{1}{16}\int_{\bar{B}_{\rho}(0)}(H-c_0)^2\dd\mu - \frac{c_0}{2} \int_{\bar{B}_{\rho}(0)} \frac{\langle h,\nu\rangle}{|h|^2}\dd\mu \\
      &\quad + \frac{1}{2\rho^2} \int_{\bar{B}_{\rho}(0)} \langle h,\nu\rangle H\dd\mu.
  \end{align}
\end{definition}
\begin{remark}\label{rem:mon-form}
    By computations in \cite[Lemma~4.1]{ruppscharrer2023} (see also \cite[Equation~(4.12)]{ruppscharrer2023}), the mapping $(0,\infty)\ni\rho\mapsto \gamma_{h,c_0}(\rho)$ is nondecreasing for any $c_0\in \R$. Moreover, for $r>0$, one computes
  \begin{equation}
      \gamma_{r h,r^{-1}c_0}(r\rho) = \gamma_{h,c_0}(\rho).
  \end{equation}
\end{remark}

With the monotonicity property and scaling behavior in \Cref{rem:mon-form}, employing arguments similar to \cite[Equation~(A.23)]{kuwertschaetzle2004}, one finds the following relation between the density at infinity and the Willmore energy of a noncompact concentration limit.

\begin{lemma}\label{lem:dens-at-infty}
    Let $A_0,V_0>0$ satisfy \eqref{eq:iso}, $c_0\in\R$, and consider a sequence of embeddings $f_j\colon\Sigma\to\R^3$ of a compact, closed and oriented surface $\Sigma$ with $\A(f_j)=A_0$, $\V(f_j)=V_0$ and
    \begin{equation}
        \sup_{j\in\N} \H_{c_0}(f_j) < \Omega(A_0,V_0,c_0).
    \end{equation}
    Suppose that, for $r_j\searrow 0$, $\hat{f}\colon\hat{\Sigma}\to\R^3$ is an immersion such that $\frac{1}{r_j}f_j\to\hat{f}$ smoothly on compact subsets of $\R^3$ after reparametrization. Then 
    \begin{equation}
        \lim_{\rho\to\infty} \frac{\hat{\mu}(\bar{B}_{\rho}(0))}{\pi\rho^2} + \frac{1}{4\pi} \W(\hat{f}) < 2.
    \end{equation}
\end{lemma}
\begin{proof}
    By smooth convergence on compact subsets of $\R^3$, one finds
    \begin{equation}
      \gamma_{\hat{f},0}(\rho) = \lim_{j\to\infty} \gamma_{\frac{1}{r_j}f_j,r_jc_1}(\rho)
    \end{equation}
    for any $c_1\in\R$. 
    Therefore, using \Cref{rem:mon-form} we have
    \begin{align}
      \lim_{\rho\to\infty} \gamma_{\hat{f},0}(\rho) &= \sup_{\rho>0} \gamma_{\hat{f},0}(\rho) = \sup_{\rho>0} \liminf_{j\to\infty }\gamma_{\frac{1}{r_j}f_j,r_jc_1}(\rho) \\
      &\leq \liminf_{j\to\infty} \sup_{\rho>0} \gamma_{\frac{1}{r_j}f_j,r_jc_1}(\rho) = \liminf_{j\to\infty} \lim_{\rho\to\infty} \gamma_{\frac{1}{r_j}f_j,r_jc_1}(\rho) \\ \label{eq:dens-at-infty}
      &= \liminf_{j\to\infty} \lim_{\rho\to\infty} \gamma_{f_j,c_1}(r_j\rho) = \liminf_{j\to\infty} \Big( \frac14\H_{c_1}(f_j) - \frac{c_1}{2}\int_{\Sigma} \frac{\langle f_j,\nu_j\rangle}{|f_j|^2}\dd\mu_j \Big).
    \end{align}
    Now, if $c_0\geq 0$ and $\Omega(A_0,V_0,c_0)= (\sqrt{8\pi}-\frac{1}{2}\sqrt{c_0^2A_0})_+^2$, we may use \eqref{eq:dens-at-infty} with $c_1=0$ and \eqref{eq:Omega_motivation2} to conclude 
    \begin{equation}\label{eq:willmore-dens-at-infty}
        \lim_{\rho\to\infty} \gamma_{\hat{f},0}(\rho) <2\pi.
    \end{equation}
    Otherwise, we use \eqref{eq:dens-at-infty} with $c_1=c_0$ and \eqref{eq:Omega_motivation0} to conclude \eqref{eq:willmore-dens-at-infty}.
    By the smooth convergence $\frac{1}{r_j}f_j\to \hat{f}$ on compact sets, one finds
    \begin{align}
      \frac14\int_{{B}_{\rho}(0)} \hat{H}^2\dd\hat{\mu} &\leq \liminf_{j\to\infty} \frac14\int_{\bar{B}_{\rho}(0)} (H_{\frac{1}{r_j}f_j}-r_jc_0)^2\dd\mu_{\frac{1}{r_j}f_j} \leq \liminf_{j\to\infty} \H_{r_jc_0}(\frac{1}{r_j}f_j) \\
      &= \liminf_{j\to\infty}\H_{c_0}(f_j) < \Omega(A_0,V_0,c_0).
    \end{align}
    Letting $\rho\to\infty$, $\W(\hat{f})<\infty$. Therefore, arguing as in \cite[Equation~(A.14)]{kuwertschaetzle2004}  
    yields
    \begin{equation}
      \lim_{\rho\to\infty} \gamma_{\hat{f},0}(\rho) = \lim_{\rho\to\infty} \frac{\hat{\mu}(\bar{B}_{\rho}(0))}{\rho^2} + \frac{1}{4}\W(\hat{f})
    \end{equation}
    and thus the claim.
\end{proof}


\begin{lemma}\label{cor:appl-of-bryant}
  Let $A_0,V_0>0$ satisfy \eqref{eq:iso}, $c_0\in\R$, and consider a sequence of spherical embeddings $f_j\colon\S^2\to\R^3$ with $\A(f_j)=A_0$, $\V(f_j)=V_0$ and
  \begin{equation}\label{eq:cor:appl-of-bryant-energy-bound}
      \sup_{j\in\N} \H_{c_0}(f_j) < \Omega(A_0,V_0,c_0).
  \end{equation}
  Suppose that $\hat{f}\colon\hat{\Sigma}\to\R^3$ is a Willmore immersion such that, 
  for some $r_j>0$ with $r_j\to r\in[0,\infty)$, $\frac{1}{r_j}f_j\to\hat{f}$ smoothly on compact subsets of $\R^3$ after reparametrization. If $\int_{B_1(0)}|\hat{A}|^2\dd\hat{\mu}>0$, then $r>0$.
\end{lemma}
\begin{proof}
    For the sake of contradiction, suppose that $r=0$. \Cref{lem:dens-at-infty} yields that
    \begin{equation}
        \lim_{\rho\to\infty} \frac{\hat{\mu}(\bar{B}_{\rho}(0))}{\pi\rho^2} + \frac{1}{4\pi} \W(\hat{f}) < 2.
    \end{equation}
    Since $\A(\frac{1}{r_j}f_j)=\frac{1}{r_j^2}A_0\to\infty$, one finds $\lim_{\rho\to\infty} \frac{\hat{\mu}(\bar{B}_{\rho}(0))}{\pi\rho^2}\geq 1$, cf. \cite[Equation~(A.22)]{kuwertschaetzle2004}. Particularly, $\W(\hat{f})<4\pi$. Moreover, by smooth convergence on compact subsets, we find that
    \begin{align}
        \int_{\hat\Sigma}|\hat{A}|^2\dd\hat{\mu}\leq \liminf_{j\to\infty} \int_{\mathbb{S}^2} |A_{f_j}|^2\dd\mu_{f_j} \leq C(A_0,V_0,c_0)<\infty,
    \end{align}
    using Gauss--Bonnet, \Cref{lem:will-vs-helf}, and \eqref{eq:cor:appl-of-bryant-energy-bound}.
    If $x_0\notin\hat{f}(\hat{\Sigma})$, consider the inversion $I(x)=|x-x_0|^{-2}(x-x_0)$ and set $\overline{\Sigma}=I(\hat{f}(\hat{\Sigma}))\cup\{0\}$. By \cite[Lemma~4.3]{kuwertschaetzle2004}, $\overline{\Sigma}$ is a smooth Willmore surface with 
    \begin{equation}
      \W(\overline{\Sigma}) = \W(\hat{f}) + 4\pi < 8\pi.
    \end{equation}
    Arguing as on \cite[p.~346]{kuwertschaetzle2004}, since all $f_j$ are spherical immersions, one finds that $\overline{\Sigma}$ necessarily also is a (topological) sphere. Thus, by \cite[Theorem~E and the Remark on p.~47]{bryant1984}, $\overline{{\Sigma}}$ is a round sphere which in turn yields that $\hat{f}(\hat{\Sigma})$ parametrizes a flat plane. This contradicts $\int_{B_1(0)}|\hat{A}|^2\dd\hat{\mu}>0$.
\end{proof}

\begin{proof}[Proof of \Cref{thm:main-result}.]
    Using \eqref{eq:en-dec}, either $\partial_tf\equiv 0$ on $[0,T)\times\S^2$ or $\H_{c_0}(f(t))<\H_{c_0}(f_0)$ for $t\in(0,T)$. Thus, replacing $f_0$ by $f(\delta)$ for some $\delta\in(0,T)$, we may without loss of generality assume \eqref{eq:life-cor-0} for some $\eta\in(0,1)$.
    \Cref{prop:concentration-limit} thus implies the existence of a concentration limit $\hat{f}\colon\hat{\Sigma}\to\R^2$. From \Cref{cor:appl-of-bryant} we infer $r>0$. Moreover, by Aleksandrov's theorem \cite{Aleksandrov1962} and Assumption \eqref{eq:iso} we have $H_{\hat f}\not\equiv const$. Hence, the conclusion follows from \Cref{lem:r>0 gives conv}.
\end{proof}

\section{Convergence for axisymmetric tori} \label{sec:tori}

For an immersion $\gamma\colon\S^1\to\mathbb{H}^2$ into the hyperbolic plane $\mathbb{H}^2=\R\times(0,\infty)$, consider the associated axisymmetric torus $f_{\gamma}\colon\mathbb{T}^2\to\R^3$ given by
\begin{equation}
    f_{\gamma}(u,v)=(\gamma^{(1)}(u),\gamma^{(2)}(u)\cos(v),\gamma^{(2)}(u)\sin(v))
\end{equation}
where $\mathbb{T}^2=\S^1\times\S^1$ with the usual identification $\S^1={\R\ }/{2\pi\mathbb{Z}}$. 


For the $(A_0,V_0,c_0)$-Helfrich flow starting in an axisymmetric torus, one has the following analog of \cite[Lemma~3.3]{dallacquamullerschatzlespener2020} with essentially the same proof with some minor modifications which we therefore postpone to \Cref{app:axi-sym}.

\begin{lemma}\label{lem:ax-sym-pres}
    Let $f_0=f_{\gamma_0}$ be an axisymmetric torus, $A_0,V_0>0$ satisfy \eqref{eq:iso} and let $c_0\in\R$. If $f\colon[0,T)\times\mathbb{T}^2\to\R^3$ is an $(A_0,V_0,c_0)$-Helfrich flow starting in $f_0$, then $f(t)=f_{\gamma(t)}$ for a smooth family of immersions $\gamma\colon[0,T)\times\S^1\to\mathbb{H}^2$.
\end{lemma}

As in \cite{dallacquamullerschatzlespener2020}, an important quantity is the so-called \emph{hyperbolic length} of $\gamma$, that is 
\begin{equation}\label{eq:def-hyp-len}
  \Ll_{\mathbb{H}^2}(\gamma)=\int_{\S^1}|\partial_u\gamma|_{g_{\mathbb{H}^2}}\dd u = \int_{\S^1} 1\dd s
\end{equation}
where $\langle\cdot,\cdot\rangle_{g_{\mathbb{H}^2}(p)}=\frac{1}{(p^{(2)})^2}\langle\cdot,\cdot\rangle_{\R^2}$ for $p=(p^{(1)},p^{(2)})\in \mathbb{H}^2$ and where $\dd s=|\partial_u\gamma|_{g_{\mathbb{H}^2}}\dd u$ is the (hyperbolic) arc-length element. In \cite{dallacquamullerschatzlespener2020}, the result in \cite{muellerspener2020} is used to control the hyperbolic length of the profile curves along the Willmore flow: For each $\varepsilon>0$, there exists $C=C(\varepsilon)>0$ such that, if $\W(f_{\gamma})\leq 8\pi-\varepsilon$, then $\Ll_{{\mathbb{H}^2}}(\gamma)\leq C(\varepsilon)$. While such a priori energy control is natural for the Willmore flow considered in \cite{dallacquamullerschatzlespener2020}, it does not immediately apply along an $(A_0,V_0,c_0)$-Helfrich flow starting in an axisymmetric torus. Therefore, in the following proposition, boundedness of the hyperbolic length of the profile curves along the evolution is shown using the Li--Yau inequality in \Cref{prp:Omega-and-LiYau}, cf. \cite[Corollary~4.11]{ruppscharrer2023}, essentially following the strategy in \cite[Proposition~4.2]{schlierf2024DirichletWillmore}.

\begin{proposition}\label{prop:bd-hyp-len}
  Let $f_0=f_{\gamma_0}$ be an embedded axisymmetric torus, let $A_0,V_0>0$ satisfy \eqref{eq:iso}, $c_0\in\R$ and suppose that
  \begin{equation}
    \H_{c_0}(f_0) < \Omega(A_0,V_0,c_0).
  \end{equation}
  If $f\colon[0,T)\times\mathbb{T}^2\to\R^3$ is an $(A_0,V_0,c_0)$-Helfrich flow starting in $f_0$ with $f(t)=f_{\gamma(t)}$ for $0\leq t<T$ and a smooth family of immersions $\gamma\colon[0,T)\times\S^1\to\mathbb{H}^2$, then
  \begin{equation}
    \sup_{0\leq t<T} \Ll_{\mathbb{H}^2}(\gamma(t)) < \infty.
  \end{equation}
\end{proposition}
\begin{proof}
  Firstly, proceeding as in \eqref{eq:diam-bd-c0geq0} and \eqref{eq:diam-bd-c0leq0}, using \Cref{rem:energy-decay} and \eqref{eq:will-vs-helf}, one finds that 
  \begin{equation}\label{eq:bd-hl-1}
    0<\gamma^{(2)}(t,\cdot)\leq \frac12\mathrm{diam}(f(t)) \leq C(A_0,c_0,\H_{c_0}(f_0)).
  \end{equation}
  Moreover, using \cite[Lemma~2.5]{dallacquamullerschatzlespener2020} with \eqref{eq:will-vs-helf} and \Cref{rem:energy-decay}, also
  \begin{equation}
    \Ll(\gamma(t)) = \int_{\S^1}|\partial_u \gamma(t,u)|\dd u \leq C(A_0,c_0,\H_{c_0}(f_0)).
  \end{equation}
  For the sake of contradiction, suppose that $\limsup_{t\nearrow T}\Ll_{\mathbb{H}^2}(\gamma(t))=\infty$. Revisiting \eqref{eq:def-hyp-len}, the above yields $\liminf_{t\nearrow T} \min_{u\in\S^1}\gamma^{(2)}(t,u) = 0$. So choose sequences $t_j\nearrow T$ and $u_j\in\S^1$ with $\gamma^{(2)}(t_j,u_j)\to 0$ and write $z_j=f(t_j,(u_j,0))\in\R^3$. Next, denote by $V_j$ the oriented varifold in $\R^3$ induced by $f_{\gamma(t_j)}$, cf. \cite[Section~2.3]{ruppscharrer2023}. Using \Cref{prp:Omega-and-LiYau,rem:energy-decay}, each $f_{\gamma(t_j)}$ is an embedding and therefore each $V_j$ belongs to the class of volume varifolds in \cite[Definition~7.4]{Scharrer23}. By \eqref{eq:a-a0-h},
  \begin{equation}
    \int_{\mathbb{T}^2}|A|^2\dd\mu = 4 \W(f(t)) \leq C(A_0,c_0,\H_{c_0}(f_0)),
  \end{equation}
  using also \eqref{eq:will-vs-helf} and \Cref{rem:energy-decay}. Due to \eqref{eq:bd-hl-1}, one can choose $h_j\in\R$ such that, for some compact $K\subseteq\R^3$, $\mathrm{supp}(V_j-(h_j,0,0))\subseteq K$ for all $j\in\N$.
  
  Thus, by \cite[Theorem~7.6]{Scharrer23}, after passing to a subsequence without relabeling, $(z_j-(h_j,0,0))\to z^*\in\R\times\{0\}^2$, $(V_j-(h_j,0,0))\to V$ as oriented varifolds for a volume varifold $V$, and 
  \begin{equation}
    \H_{c_0}(V)\leq \liminf_{j\to\infty}\H_{c_0}(V_j) < \Omega(A_0,V_0,c_0).
  \end{equation}
 
  Particularly, using 
  \cite[Corollary~4.11]{ruppscharrer2023}, 
  \begin{equation}\label{eq:bd-hl-2}
    \theta^2(\mu_V,z) < 2\quad\text{for all $z\in\R^3$}.
  \end{equation}
  Consider now reparametrizations $\widetilde{\gamma}_j\colon\S^1\to\mathbb{H}^2$ of $\gamma(t_j)-(h_j,0)$ by constant Euclidean speed, that is, such that $\partial_u|\partial_u\widetilde{\gamma}_j|=0$ on $\S^1$. After passing to a further subsequence without relabeling, using \cite[Lemma~3.7]{schlierf2024DirichletWillmore}, 
  \begin{align}
    &\widetilde{\gamma}_j \rightharpoonup^* \widetilde{\gamma}\text{ in }W^{1,\infty}(\S^1) \text{ and } \widetilde{\gamma}_j \to \widetilde{\gamma}\quad\text{uniformly};\\
    &\widetilde{\gamma}_j \rightharpoonup \widetilde{\gamma}\quad\text{in $W^{2,2}([a,b])$ for every $[a,b]\subseteq\{\widetilde{\gamma}^{(2)}>0\}$}
  \end{align}
  for some $\widetilde{\gamma}\in W^{1,\infty}(\S^1,\R\times[0,\infty))$. Now proceeding as in \cite[proof of Proposition~4.2]{schlierf2024DirichletWillmore}, one argues that, for $z^*=\lim_{j\to\infty}(z_j-(h_j,0,0))\in\R\times\{0\}^2$,
  \begin{equation}
    \theta^2(\mu_V,z^*)\geq 2,
  \end{equation}
  contradicting \eqref{eq:bd-hl-2}.
\end{proof}


\begin{proof}[Proof of \Cref{intro:cor:axisymmetric_tori}.]
  As for \Cref{thm:main-result}, one may without loss of generality assume \eqref{eq:life-cor-0} for some $\eta\in (0,1)$. Then consider a concentration limit $\hat{f}\colon\hat{\Sigma}\to\R^3$ as in \Cref{prop:concentration-limit} and proceed as in \cite[proof of Theorem~1.2]{dallacquamullerschatzlespener2020} to deduce from \Cref{prop:bd-hyp-len} that $\hat{\Sigma}$ is necessarily compact and $r>0$. Clearly $H_{\hat{f}}\not\equiv const$ by Aleksandrov's theorem \cite{Aleksandrov1962}. So the claim follows from \Cref{lem:r>0 gives conv}, arguing as in \cite[Lemma~3.4]{dallacquamullerschatzlespener2020}.
\end{proof}

\appendix 

\section{Proof of \texorpdfstring{\Cref{lem:loc-2}}{Lemma 3.1}}\label{app:local}

This section is devoted to proving \Cref{lem:loc-2}. We first compute the exact localized evolution of the energy.

\begin{lemma}\label{lem:loc-1}
  Let $A_0,V_0>0$ satisfy \eqref{eq:iso}, $c_0\in\R$ and $f\colon[0,T)\times\Sigma\to\R^3$ be an $(A_0,V_0,c_0)$-Helfrich flow. If $\widetilde{\eta}\in C_c^{\infty}(\R^3)$ and $\eta=\widetilde{\eta}\circ f$, 
  \begin{align}
    \frac{\dd}{\dd t}&\int_{\Sigma}\frac12 H^2\eta\dd \mu + \int_{\Sigma} |\nabla\W_0(f)|^2\eta\dd\mu\\
    &= (\lambda_1+\frac12c_0^2)\int_{\Sigma}\Delta H H\eta\dd\mu + \int_{\Sigma}\big((\lambda_1+\frac12c_0^2)H+\lambda_2\big)|A^0|^2H\eta\dd \mu\\
    &\quad-2\int_{\Sigma}\langle\nabla\W_0(f),\nu\rangle\langle\nabla H,\nabla\eta\rangle_{g_f}\dd\mu - \int_{\Sigma}\langle\nabla\W_0(f),\nu\rangle H\Delta\eta\dd\mu + \frac12\int_{\Sigma}H^2\partial_t\eta\dd\mu\\
    &\quad+c_0\int_{\Sigma}(|A^0|^2-\frac12H^2)(\langle\nabla\W_0(f),\nu\rangle\eta + 2\langle \nabla H,\nabla\eta\rangle_{g_f}+H\Delta\eta)\dd\mu
  \end{align}
  as well as  
  \begin{align}
    \partial_t&\int_{\Sigma}|A^0|^2\eta\dd\mu+\int_{\Sigma}|\nabla\W_0(f)|^2\eta\dd\mu\\
    &= (2\lambda_1+c_0^2)\int_{\Sigma}\langle\nabla^2H,A^0\rangle_{g_f}\eta\dd\mu + \int_{\Sigma}\big((\lambda_1+\frac12c_0^2)H+\lambda_2\big)|A^0|^2H\eta\dd\mu\\
    &\quad-2\int_{\Sigma}\langle\nabla\W_0(f),\nu\rangle\bigl(\langle A^0,\nabla^2\eta\rangle_{g_f}+\langle \nabla H,\nabla\eta\rangle_{g_f}\bigr)\dd\mu \\
    &\quad + c_0\int_{\Sigma}(|A^0|^2-\frac12H^2) \bigl( \langle\nabla\W_0(f),\nu\rangle\eta + 2\langle A^0,\nabla^2\eta\rangle_{g_f}+2\langle \nabla H,\nabla\eta\rangle_{g_f} \bigr)\dd\mu\\
    &\quad + \int_{\Sigma}|A^0|^2\partial_t\eta\dd\mu.
  \end{align}
\end{lemma}
\begin{proof}
  Using \eqref{eq:ev-mu}, \eqref{eq:ev-H} and \cite[(31) and (32)]{kuwertschaetzle2001}, writing $\partial_tf=\xi\nu$,
  \begin{align}
    \partial_t&\int_{\Sigma}\frac12H^2\eta\dd\mu + \int_{\Sigma}|\nabla\W_0(f)|^2\eta\dd\mu\\
    &= \int_{\Sigma}\big((\lambda_1+\frac12c_0^2)H+\lambda_2\big)(\Delta H+|A^0|^2H)\eta\dd\mu+\int_{\Sigma} (2\xi\langle\nabla H,\nabla\eta\rangle_{g_f}+H\xi\Delta\eta)\dd\mu\\
    &\quad+c_0\int_{\Sigma}(|A^0|^2-\frac12H^2)(\Delta H+|A^0|^2H)\eta\dd\mu+\frac12\int_{\Sigma}H^2\partial_t\eta\dd\mu,
  \end{align}
  since $\xi=-(\Delta H+|A^0|^2H)+c_0(|A^0|^2-\frac12H^2)+(\lambda_1+\frac12c_0^2)H+\lambda_2$. Moreover, also using $\Delta(H\eta) = \Delta H\eta + 2\langle \nabla H,\nabla \eta\rangle_{g_f}+H\Delta \eta$, one finds
  \begin{align}
    \partial_t&\int_{\Sigma}\frac12H^2\eta\dd\mu + \int_{\Sigma}|\nabla\W_0(f)|^2\eta\dd\mu\\
    &= \int_{\Sigma}\big((\lambda_1+\frac12c_0^2)H+\lambda_2\big)(\Delta (H\eta)+|A^0|^2H\eta)\dd\mu\\
    &\quad-2\int_{\Sigma}\langle\nabla\W_0(f),\nu\rangle\langle\nabla H,\nabla\eta\rangle_{g_f}\dd\mu-\int_{\Sigma}\langle\nabla\W_0(f),\nu\rangle H\Delta\eta\dd\mu\\
    &\quad+c_0\int_{\Sigma}(|A^0|^2-\frac12H^2)(\langle\nabla\W_0(f),\nu\rangle\eta + 2\langle \nabla H,\nabla\eta\rangle_{g_f}+H\Delta\eta)\dd\mu+\frac12\int_{\Sigma}H^2\partial_t\eta\dd\mu.
  \end{align}
  The first claim follows after integrating by parts in the first term on the right. Moreover, as in 
  \cite[Article~C, proof of Lemma~2.8]{rupp2022PhD}
  and \cite[p.~423]{kuwertschaetzle2001}, one finds
  \begin{align}
    \partial_t(|A^0|^2\dd\mu) &= 2\nabla_i(\nabla_j\xi A^0(e_i,e_j))\dd\mu-\nabla_j\xi\nabla_jH\dd\mu+|A^0|^2H\xi\dd\mu\\
    &=2\nabla_i(\nabla_j\xi A^0(e_i,e_j))\dd\mu-\nabla_j(\xi\nabla_jH)\dd\mu+(\Delta H+|A^0|^2H)\xi\dd\mu.
  \end{align}
  Again, with $\xi=-(\Delta H+|A^0|^2H)+c_0(|A^0|^2-\frac12H^2)+(\lambda_1+\frac12c_0^2)H+\lambda_2$, 
  \begin{align}
    \partial_t&(|A^0|^2\dd\mu) + |\nabla\W_0(f)|^2\dd\mu \\
    &= 2\nabla_i(\nabla_j\xi A^0(e_i,e_j))\dd\mu-\nabla_j(\xi\nabla_jH)\dd\mu\\
    &\quad + c_0(|A^0|^2-\frac12H^2)\langle\nabla\W_0(f),\nu\rangle \dd\mu+\big((\lambda_1+\frac12c_0^2)H+\lambda_2\big)\langle\nabla\W_0(f),\nu\rangle\dd\mu.
  \end{align}
  Integrating by parts and using $\nabla_iH=2(\nabla_jA^0)_{ij}$ by Codazzi--Mainardi, one obtains
  \begin{align}
    \partial_t&\int_{\Sigma}|A^0|^2\eta\dd\mu+\int_{\Sigma}|\nabla\W_0(f)|^2\eta\dd\mu\\
    &= \int_{\Sigma} 2\xi A^0_{ij}\nabla_{ij}^2\eta+2\xi\nabla_jH\nabla_j\eta + \big((\lambda_1+\frac12c_0^2)H+\lambda_2\big)\langle\nabla\W_0(f),\nu\rangle\eta\dd\mu\\
    &\quad +\int_{\Sigma}c_0(|A^0|^2-\frac12H^2)\langle\nabla\W_0(f),\nu\rangle \eta \dd\mu + \int_{\Sigma} |A^0|^2\partial_t\eta\dd\mu.
  \end{align}
  Integration by parts yields
  \begin{equation}
    \int_{\Sigma}\big((\lambda_1+\frac12c_0^2)H+\lambda_2\big)(2A^0_{ij}\nabla^2_{ij}\eta+2\nabla_jH\nabla_j\eta+\Delta H\eta)\dd\mu = (2\lambda_1+c_0^2)\int_{\Sigma}\langle\nabla^2H,A^0\rangle_{g_f}\eta\dd\mu,
  \end{equation}
  so that the claim follows.
\end{proof}

\begin{proof}[{Proof of \Cref{lem:loc-2}}]
  Using \eqref{eq:a-a0-h} and \Cref{lem:loc-1} with $\tilde{\eta}\vcentcolon = \tilde{\gamma}^4$, so that $\eta = \gamma^4$, one finds
  \begin{align}
    &\frac{\dd}{\dd t}\int_{\Sigma}|A|^2\gamma^4\dd\mu+2\int_{\Sigma}|\nabla\W_0(f)|^2\gamma^4\dd\mu\\
    &= (2\lambda_1+c_0^2)\int_{\Sigma}\bigl(\langle\nabla^2H,A\rangle_{g_f}+|A^0|^2H^2\bigr)\gamma^4\dd\mu+2\lambda_2 \int_{\Sigma}|A^0|^2H\gamma^4\dd\mu\\
    &\quad -2\int_{\Sigma}\langle\nabla\W_0(f),\nu\rangle\bigl(2\langle\nabla H,\nabla\gamma^4\rangle_{g_f}+\langle\nabla^2\gamma^4,A\rangle_{g_f}\bigr)\dd\mu+\int_{\Sigma}|A|^2\partial_t\gamma^4\dd\mu\\
    &\quad + c_0\int_{\Sigma}(|A^0|^2-\frac12H^2)\bigl( 2\langle\nabla^2\gamma^4,A\rangle_{g_f}+4\langle\nabla H,\nabla\gamma^4\rangle_{g_f}+2\langle\nabla\W_0(f),\nu\rangle\gamma^4 \bigr)\dd\mu.
  \end{align}
  We now estimate the individual terms starting with $\int_\Sigma|A|^2\partial_t\gamma^4\dd\mu$. A straight-forward computation yields 
  \begin{equation}
    |\partial_t\gamma^4|\leq C\Lambda\gamma^3\bigl( |\nabla\W_0(f)|+(|\lambda_1|+c_0^2)|A|+|\lambda_2|+|c_0||A|^2 \bigr).
  \end{equation}
  The last term arising in $\int_\Sigma|A|^2 \partial_t \gamma^4 \dd\mu$ is estimated by
  \begin{equation}
      C\Lambda|c_0|\int_\Sigma|A|^4\gamma^3\dd\mu\le C c_0^2\int_\Sigma|A|^4\gamma^4\dd\mu + C\Lambda^2\int_\Sigma|A|^4\gamma^2\dd\mu
  \end{equation}
  whereas the first term can be absorbed via 
  \begin{equation}
      C\Lambda\int_\Sigma|A|^2|\nabla\mathcal W_0(f)|\gamma^3\dd \mu \le \varepsilon\int_\Sigma|\nabla\mathcal W_0(f)|^2\gamma^4\dd\mu +C(\varepsilon)\Lambda^2\int_{\Sigma}|A|^4\gamma^2\dd\mu.
  \end{equation}
  Moreover, as  $|\nabla^2\gamma^4|\leq C(\Lambda^2+|A|\Lambda\gamma)\gamma^2$,
  \begin{align}
    \int_{\Sigma} |\nabla\W_0(f)||\nabla^2\gamma^4||A|\dd\mu &\leq \varepsilon \int_{\Sigma}|\nabla\W_0(f)|^2\gamma^4\dd\mu + C(\varepsilon)\Lambda^2\int_{\Sigma}|A|^4\gamma^2\dd\mu\\
    &\quad+C(\varepsilon)\Lambda^4\int_{\{\gamma>0\}}|A|^2\dd\mu.
  \end{align}
  Furthermore, as in \cite[proof of Lemma 3.2]{kuwertschaetzle2001}, one obtains
  \begin{align}
    \int_{\Sigma}|\nabla H|^2\gamma^2\dd\mu&\leq \varepsilon \int_{\Sigma} |\nabla\W_0(f)|^2\gamma^4\dd\mu \\
    &\quad+ C\Big(\Lambda^2+\frac{1}{\varepsilon}\Big)\int_{\{\gamma>0\}} H^2\dd\mu + C\int_{\Sigma}|A|^4\gamma^2\dd\mu.\label{eq:loc-2-1}
  \end{align}
  Therefore, using $|\nabla\gamma^4|\leq C\gamma^3\Lambda$,
  \begin{align}
     \int_{\Sigma}|\nabla\W_0(f)||\nabla H||\nabla \gamma^4|\dd\mu &\leq \varepsilon \int_{\Sigma}|\nabla\W_0(f)|^2\gamma^4\dd\mu+C(\varepsilon)\Lambda^2\int_{\Sigma}|A|^4\gamma^2\dd\mu\\
     &\quad+C(\varepsilon)\Lambda^4\int_{\{\gamma>0\}}H^2\dd\mu.
  \end{align}
  Using $|c_0(|A^0|^2-\frac12H^2)|\leq C |c_0||A|^2$ as well as $|\nabla^2\gamma^4|\leq C(\Lambda^2+|A|\Lambda\gamma)\gamma^{2}$, the terms containing $c_0$ can be dealt with as follows. Firstly, we have
  \begin{align}
      |c_0|\int_\Sigma|A|^3|\nabla^2\gamma^4|\dd\mu & \le C|c_0|\Lambda^2\int_\Sigma|A|^3\gamma^2\dd\mu + C|c_0|\Lambda\int_\Sigma|A|^4\gamma^3\dd\mu \\
      &\le C c_0^2\int_\Sigma|A|^4\gamma^4\dd\mu + C\Lambda^4\int_{\{\gamma>0\}}|A|^2\dd\mu + C\Lambda^2\int_\Sigma|A|^4\gamma^2\dd\mu.
  \end{align}
  Secondly, using \eqref{eq:loc-2-1},
  \begin{align}
      |c_0|\int_\Sigma|A|^2|\nabla H||\nabla\gamma^4|\dd \mu& \le Cc_0^2\int_\Sigma|A|^4\gamma^4\dd\mu +\varepsilon\int_\Sigma|\nabla\mathcal W_0(f)|^2\gamma^4 \dd\mu\\
      &\quad + C(\varepsilon)\Lambda^2\int_\Sigma|A|^4\gamma^2\dd\mu + C(\varepsilon)\Lambda^4\int_{\{\gamma>0\}}|A|^2\dd\mu 
  \end{align}
  and finally,
  \begin{equation}
      |c_0|\int_\Sigma |A|^2|\nabla \mathcal W_0(f)|\gamma^4\dd\mu \le \varepsilon\int_\Sigma|\nabla\mathcal W_0(f)|^2\gamma^4\dd\mu + C(\varepsilon)c_0^2\int_\Sigma|A|^4\gamma^4\dd\mu
  \end{equation}
  which completes the proof after taking $\varepsilon>0$ sufficiently small and absorbing.
\end{proof}

\section{On preservation of axisymmetry}\label{app:axi-sym}

Following \cite[Definition~1.1]{dallacquamullerschatzlespener2020} with the identification $\S^1=\R\ /2\pi\mathbb{Z}$, a smooth immersion $f\colon\mathbb{T}^2\to\R^3$ is said to be an axisymmetric torus if, for a smooth immersion $\gamma\colon\S^1\to\mathbb{H}^2$,
\begin{equation}
    f(u,v)=(\gamma^{(1)}(u),\gamma^{(2)}(u)\cos(v),\gamma^{(2)}(u)\sin(v))
\end{equation}
where $\mathbb{T}^2=\S^1\times\S^1$. One has the following characterization.

\begin{proposition}[{\cite[Proposition~3.2]{dallacquamullerschatzlespener2020}}]\label{prop:3.2indmss}
    A smooth immersion $f\colon\mathbb{T}^2\to\R^3$ is an axisymmetric torus if and only if
    \begin{enumerate}[(a)]
        \item\label{item:char-axi-a} $f(u,v+\varphi) = R_{\varphi}f(u,v)$ for all $u,v,\varphi\in\S^1$ where
        \begin{equation}
            R_{\varphi} = \begin{pmatrix}
                1&0&0\\
                0&\cos\varphi&-\sin\varphi\\
                0&\sin\varphi&\cos\varphi
            \end{pmatrix}\quad\text{and}
        \end{equation}
        \item\label{item:char-axi-b} $f^{(3)}(u,0)=0$ for all $u\in\S^1$ and $f^{(2)}(u_0,0)\geq 0$ for at least one $u_0\in\S^1$.
    \end{enumerate}
\end{proposition}

For \Cref{lem:ax-sym-pres}, if one followed the proof of \cite[Lemma~3.3]{dallacquamullerschatzlespener2020} line-by-line, in the case of the Helfrich flow, one would have to use \Cref{cor:life,prop:ho-en-est} in place of \cite[Theorem~1.2 and Equation~(4.27)]{kuwertschaetzle2002}. Then however, one needs to additionally assume \eqref{eq:en-thresh}. Therefore, we repeat the general strategy of the proof of \cite[Lemma~3.3]{dallacquamullerschatzlespener2020} but avoid using \Cref{cor:life,prop:ho-en-est}, i.e., do not assume \eqref{eq:en-thresh}.

\begin{proof}[Proof of \Cref{lem:ax-sym-pres}.]
    Let $f\colon[0,T)\times\mathbb{T}^2\to\R^3$ be an $(A_0,V_0,c_0)$-Helfrich flow starting in an axisymmetric torus $f_0=f_{\gamma_0}$ as in the statement. We first verify \Cref{prop:3.2indmss}\eqref{item:char-axi-a} along the evolution. To this end, for fixed $\varphi\in\S^1$, $R_{\varphi}\colon\R^3\to\R^3$ is an isometry. Therefore, one readily checks that $(R_{\varphi}^{-1})f(t,\cdot+(0,\varphi))\colon[0,T)\times\mathbb{T}^2\to\R^3$ is an $(A_0,V_0,c_0)$-Helfrich flow starting in $R_{\varphi}^{-1}f_0(\cdot+(0,\varphi))=f_0$, using that $f_0$ is an axisymmetric torus. Well-posedness of \eqref{eq:flow-eq}, see \cite{KohsakaNagasawa2006}, then yields that each $f(t)$ satisfies \Cref{prop:3.2indmss}\eqref{item:char-axi-a} for $0\leq t<T$.

    Therefore, there exist smooth functions $x,y,z\colon[0,T)\times\S^1\to\R$ with
    \begin{equation}\label{eq:ax-sym-pres-1}
        f(t,(u,v)) = R_vf(t,(u,0)) = R_v \big(x(t,u),y(t,u),z(t,u)\big).
    \end{equation}
    In order to establish \Cref{prop:3.2indmss}\eqref{item:char-axi-b}, we first show $z\equiv 0$ by slightly modifying the proof of \cite[Lemma~3.3]{dallacquamullerschatzlespener2020}. Consider the normal field
    \begin{equation}\label{eq:ax-sym-pres-2}
        \widetilde{\nu} \vcentcolon= \frac{\partial_uf\times\partial_vf}{|\partial_uf\times\partial_vf|} = \frac{1}{\sqrt{\det (Df^tDf)}} R_v\begin{pmatrix}
            y\partial_uy+z\partial_uz\\-y\partial_ux\\-z\partial_ux
        \end{pmatrix}.
    \end{equation}
    For the sake of contradiction, suppose that
    \begin{equation}\label{eq:ax-sym-pres-3}
        S=\sup\{s\in[0,T):f(t)\text{ is an axisymmetric torus for all }0\leq t\leq s\}<T.
    \end{equation}
    As in the proof of \cite[Lemma~3.3]{dallacquamullerschatzlespener2020}, $f(S)$ is an axisymmetric torus and $z(S,\cdot)=0$. Moreover, since $\partial_tf = \langle\partial_tf,\nu\rangle\nu$, \eqref{eq:ax-sym-pres-1} and \eqref{eq:ax-sym-pres-2} yield that on $[0,T)\times \S^1$ we have
    \begin{equation}
        \partial_t z^2 = 2z\partial_tz=\left.2z\langle \partial_t f,\widetilde{\nu}\rangle \widetilde{\nu}^{(3)}\right\vert_{v=0} = \left.2z \langle\partial_tf,\widetilde{\nu}\rangle\frac{1}{\sqrt{\det (Df^tDf)}}\right\vert_{v=0} (-z\partial_ux).
    \end{equation}
    For $0<\varepsilon<T-S$, since $f$ is a \emph{smooth} family of \emph{immersions} on $[S,S+\varepsilon]\times\mathbb{T}^2$, we find
    \begin{equation}
        \partial_t z^2 \leq \sup_{t\in[S,S+\varepsilon],u\in\S^1} \Big|\frac{2\langle\partial_tf(t,(u,0)),\widetilde{\nu}(t,(u,0))\rangle\partial_ux(t,u)}{\sqrt{\det [Df^t(t,(u,0))Df(t,(u,0))]}}\Big| \cdot z^2=\vcentcolon C\cdot z^2\quad\text{on $[S,S+\varepsilon]\times\mathbb{S}^1$}.
    \end{equation}
    Since $z(S,\cdot)\equiv 0$, this yields $z(t,u)=0$ for all $S\leq t\leq S+\varepsilon$ and $u\in\S^1$. 
    Notice that $y(S,\cdot)>0$ since $f(S)$ is an axisymmetric torus. Therefore, for $\varepsilon>0$ small enough, there exists $u_0\in \mathbb S^1$ such that $y(t,u_0)>0$ for all $S\le t \le S + \varepsilon$. Thus, by \Cref{prop:3.2indmss}, $f(t)$ is an axisymmetric torus for $t\in [0,S+\varepsilon]$, contradicting \eqref{eq:ax-sym-pres-3}.
\end{proof}

\section*{Acknowledgments}
This research was funded in whole, or in part, by the Austrian Science Fund (FWF), grant numbers \href{https://doi.org/10.55776/P32788}{10.55776/P32788} and \href{https://doi.org/10.55776/ESP557}{10.55776/ESP557}. 

\section*{Data availability} 
Data sharing not applicable to this article as no datasets were generated or analyzed during the current study. 

\bibliographystyle{alpha}
\bibliography{biblio}

@article {KubinLussardiMorandotti24,
    AUTHOR = {Kubin, Anna and Lussardi, Luca and Morandotti, Marco},
     TITLE = {Direct minimization of the {C}anham-{H}elfrich energy on
              generalized {G}auss graphs},
   JOURNAL = {J. Geom. Anal.},
  FJOURNAL = {Journal of Geometric Analysis},
    VOLUME = {34},
      YEAR = {2024},
    NUMBER = {5},
     PAGES = {Paper No. 121, 25},
      ISSN = {1050-6926,1559-002X},
   MRCLASS = {49Q20 (49J45 49Q10 53C80 92C10)},
  MRNUMBER = {4715387},
       DOI = {10.1007/s12220-024-01564-2},
       URL = {https://doi.org/10.1007/s12220-024-01564-2},
}

@article {BrazdaLussardiStefanelli20,
    AUTHOR = {Brazda, Katharina and Lussardi, Luca and Stefanelli, Ulisse},
     TITLE = {Existence of varifold minimizers for the multiphase
              {C}anham-{H}elfrich functional},
   JOURNAL = {Calc. Var. Partial Differential Equations},
  FJOURNAL = {Calculus of Variations and Partial Differential Equations},
    VOLUME = {59},
      YEAR = {2020},
    NUMBER = {3},
     PAGES = {Paper No. 93, 26},
      ISSN = {0944-2669,1432-0835},
   MRCLASS = {49Q20 (49J45 49Q10 53C80 92C10)},
  MRNUMBER = {4098040},
MRREVIEWER = {Chun-Chi\ Lin},
       DOI = {10.1007/s00526-020-01759-9},
       URL = {https://doi.org/10.1007/s00526-020-01759-9},
}

@article {KusnerMcGrath23,
    AUTHOR = {Kusner, Robert and McGrath, Peter},
     TITLE = {On the {C}anham problem: bending energy minimizers for any
              genus and isoperimetric ratio},
   JOURNAL = {Arch. Ration. Mech. Anal.},
  FJOURNAL = {Archive for Rational Mechanics and Analysis},
    VOLUME = {247},
      YEAR = {2023},
    NUMBER = {1},
     PAGES = {Paper No. 10, 14},
      ISSN = {0003-9527,1432-0673},
   MRCLASS = {49Q10 (58E12)},
  MRNUMBER = {4537344},
MRREVIEWER = {Alexis\ Michelat},
       DOI = {10.1007/s00205-022-01833-w},
       URL = {https://doi.org/10.1007/s00205-022-01833-w},
}

@article {MondinoScharrer2023,
    AUTHOR = {Mondino, Andrea and Scharrer, Christian},
     TITLE = {A strict inequality for the minimization of the {W}illmore
              functional under isoperimetric constraint},
   JOURNAL = {Adv. Calc. Var.},
  FJOURNAL = {Advances in Calculus of Variations},
    VOLUME = {16},
      YEAR = {2023},
    NUMBER = {3},
     PAGES = {529--540},
      ISSN = {1864-8258,1864-8266},
   MRCLASS = {49Q10 (53A05)},
  MRNUMBER = {4609797},
       DOI = {10.1515/acv-2021-0002},
       URL = {https://doi.org/10.1515/acv-2021-0002},
}

@article {KellerMondinoRiviere2014,
    AUTHOR = {Keller, Laura Gioia Andrea and Mondino, Andrea and Rivi\`ere,
              Tristan},
     TITLE = {Embedded surfaces of arbitrary genus minimizing the {W}illmore
              energy under isoperimetric constraint},
   JOURNAL = {Arch. Ration. Mech. Anal.},
  FJOURNAL = {Archive for Rational Mechanics and Analysis},
    VOLUME = {212},
      YEAR = {2014},
    NUMBER = {2},
     PAGES = {645--682},
      ISSN = {0003-9527,1432-0673},
   MRCLASS = {49Q05 (53C42)},
  MRNUMBER = {3176354},
MRREVIEWER = {Christine\ Breiner},
       DOI = {10.1007/s00205-013-0694-9},
       URL = {https://doi.org/10.1007/s00205-013-0694-9},
}

@article {KohsakaNagasawa2006,
    AUTHOR = {Kohsaka, Yoshihito and Nagasawa, Takeyuki},
     TITLE = {On the existence of solutions of the {H}elfrich flow and its
              center manifold near spheres},
   JOURNAL = {Differential Integral Equations},
  FJOURNAL = {Differential and Integral Equations. An International Journal
              for Theory \& Applications},
    VOLUME = {19},
      YEAR = {2006},
    NUMBER = {2},
     PAGES = {121--142},
      ISSN = {0893-4983},
   MRCLASS = {58E12 (49Q10 74G05 74L15 92C50)},
  MRNUMBER = {2194500},
MRREVIEWER = {Marc\ Michel\ Soret},
}

@article {nagasawayi2012,
    AUTHOR = {Nagasawa, Takeyuki and Yi, Taekyung},
     TITLE = {Local existence and uniqueness for the {$n$}-dimensional
              {H}elfrich flow as a projected gradient flow},
   JOURNAL = {Hokkaido Math. J.},
  FJOURNAL = {Hokkaido Mathematical Journal},
    VOLUME = {41},
      YEAR = {2012},
    NUMBER = {2},
     PAGES = {209--226},
      ISSN = {0385-4035},
   MRCLASS = {53Cxx (35K30 49Q10)},
  MRNUMBER = {2977045},
       DOI = {10.14492/hokmj/1340714413},
       URL = {https://doi.org/10.14492/hokmj/1340714413},
}

@article {schlierf2024helfrich,
    AUTHOR = {Schlierf, Manuel},
     TITLE = {Spontaneous curvature effects of the {H}elfrich flow:
              singularities and convergence},
   JOURNAL = {Comm. Partial Differential Equations},
  FJOURNAL = {Communications in Partial Differential Equations},
    VOLUME = {50},
      YEAR = {2025},
    NUMBER = {3},
     PAGES = {441--476},
      ISSN = {0360-5302,1532-4133},
   MRCLASS = {53E40 (35B40 35K41 49Q10)},
  MRNUMBER = {4870993},
       DOI = {10.1080/03605302.2025.2457047},
       URL = {https://doi.org/10.1080/03605302.2025.2457047},
}

@article {julin2024sharpquantitativealexandrovinequality,
    AUTHOR = {Julin, Vesa and Morini, Massimiliano and Oronzio, Francesca
              and Spadaro, Emanuele},
     TITLE = {A sharp quantitative {A}lexandrov inequality and applications
              to volume preserving geometric flows in 3{D}},
   JOURNAL = {Arch. Ration. Mech. Anal.},
  FJOURNAL = {Archive for Rational Mechanics and Analysis},
    VOLUME = {249},
      YEAR = {2025},
    NUMBER = {6},
     PAGES = {Paper No. 78, 45},
      ISSN = {0003-9527,1432-0673},
   MRCLASS = {35R35 (35B40 53E10)},
  MRNUMBER = {4991147},
MRREVIEWER = {Ao\ Sun},
       DOI = {10.1007/s00205-025-02141-9},
       URL = {https://doi.org/10.1007/s00205-025-02141-9},
}

@article {Ciraolo2021,
    AUTHOR = {Ciraolo, Giulio},
     TITLE = {Quantitative estimates for almost constant mean curvature
              hypersurfaces},
   JOURNAL = {Boll. Unione Mat. Ital.},
  FJOURNAL = {Bollettino dell'Unione Matematica Italiana},
    VOLUME = {14},
      YEAR = {2021},
    NUMBER = {1},
     PAGES = {137--150},
      ISSN = {1972-6724,2198-2759},
   MRCLASS = {53A10 (35B35 35B50 35B51 35J70 35N25 53C24)},
  MRNUMBER = {4218626},
       DOI = {10.1007/s40574-020-00242-9},
       URL = {https://doi.org/10.1007/s40574-020-00242-9},
}

@article{metsch2023axially,
  title={Axially Symmetric {W}illmore Minimizers with Prescribed Isoperimetric Ratio},
  author={Metsch, Jan-Henrik},
  year={2025},
  journal={Partial Differ. Equ. Appl.},
  volume = {6},
  pages      = {Paper No. 5},
  Number = {1},
  MRNUMBER = {4859183},
  DOI = {10.1007/s42985-024-00303-0},
}

@book {AndrewsChowGuentherLangford,
	AUTHOR = {Andrews, Ben and Chow, Bennett and Guenther, Christine and
	Langford, Mat},
	TITLE = {Extrinsic geometric flows},
	SERIES = {Graduate Studies in Mathematics},
	VOLUME = {206},
	PUBLISHER = {American Mathematical Society, Providence, RI},
	YEAR = {\copyright 2020},
	PAGES = {xxviii+759},
	ISBN = {978-1-4704-5596-5},
	MRCLASS = {53Exx (35K20 52A20 53A07 58J35)},
	MRNUMBER = {4249616},
	MRREVIEWER = {Greg\'{o}rio Manoel Silva Neto},
	DOI = {10.1090/gsm/206},
	URL = {https://doi.org/10.1090/gsm/206},
}

@article{schlierf2024DirichletWillmore,
  url         = {https://doi.org/10.1515/acv-2024-0018},
  title       = {Global existence for the {W}illmore flow with boundary via {S}imon’s {L}i–{Y}au inequality},
  author      = {Schlierf, Manuel},
  journal     = {Adv. Calc. Var.},
  doi         = {10.1515/acv-2024-0018},
  year        = {2025},
  lastchecked = {2025-02-14}
}

@article {RuppScharrer24,
    AUTHOR = {Rupp, Fabian and Scharrer, Christian},
     TITLE = {Global regularity of integral 2-varifolds with square
              integrable mean curvature},
   JOURNAL = {J. Math. Pures Appl. (9)},
  FJOURNAL = {Journal de Math\'ematiques Pures et Appliqu\'ees. Neuvi\`eme
              S\'erie},
    VOLUME = {204},
      YEAR = {2025},
     PAGES = {Paper No. 103797, 30},
      ISSN = {0021-7824,1776-3371},
   MRCLASS = {49Q20 (53A05 53C42)},
  MRNUMBER = {4962374},
       DOI = {10.1016/j.matpur.2025.103797},
       URL = {https://doi.org/10.1016/j.matpur.2025.103797},
}

@article{Scharrer23,
  title={Properties of surfaces with spontaneous curvature},
  author={Scharrer, Christian},
  journal={arXiv:2310.19935 [math.DG]},
  eprint        = {2310.19935},
  archiveprefix = {arXiv},
  year={2023},
  Note = {To appear in Anal. PDE},
}

@phdthesis{rupp2022PhD,
  title={Constrained gradient flows for Willmore-type functionals},
  author={Rupp, Fabian},
  year={2022},
  school={Universit{\"a}t Ulm}
}

@article{barrettgarckenuernberg2021,
  author   = {Barrett, John W. and Garcke, Harald and N\"{u}rnberg, Robert},
  title    = {Stable approximations for axisymmetric {W}illmore flow for
              closed and open surfaces},
  journal  = {ESAIM Math. Model. Numer. Anal.},
  fjournal = {ESAIM. Mathematical Modelling and Numerical Analysis},
  volume   = {55},
  year     = {2021},
  number   = {3},
  pages    = {833--885},
  issn     = {2822-7840,2804-7214},
  mrclass  = {65M60 (35K55 53E10 65M12)},
  mrnumber = {4253162},
  doi      = {10.1051/m2an/2021014},
  url      = {https://doi.org/10.1051/m2an/2021014}
}

@book{Blaschke1929,
  author    = {Wilhelm Blaschke},
  title={Vorlesungen {\"u}ber Differentialgeometrie und geometrische Grundlagen von Einsteins Relativit{\"a}tstheorie III: Differentialgeometrie der Kreise und Kugeln},
  publisher = {Springer-Verlag},
  address   = {Berlin},
  year      = {1929},
  jfm       = {55.0422.01}
}

@article{Thomsen1923,
  author  = {Thomsen, G.},
  title   = {{\"U}ber konforme {G}eometrie {I}: Grundlagen der konformen {F}l{\"a}chentheorie},
  journal = {Hamb. Math. Abh.},
  volume  = {3},
  pages   = {31--56},
  year    = {1923},
  doi     = {10.1007/BF02954615},
  jfm     = {49.0530.02}
}

@article {Jakob2023,
    AUTHOR = {Jakob, Ruben},
     TITLE = {The {W}illmore flow of {H}opf-tori in the 3-sphere},
   JOURNAL = {J. Evol. Equ.},
  FJOURNAL = {Journal of Evolution Equations},
    VOLUME = {23},
      YEAR = {2023},
    NUMBER = {4},
     PAGES = {Paper No. 72, 42},
      ISSN = {1424-3199,1424-3202},
   MRCLASS = {53E40 (35R01 53C42 58J35)},
  MRNUMBER = {4660741},
       DOI = {10.1007/s00028-023-00923-w},
       URL = {https://doi.org/10.1007/s00028-023-00923-w},
}

@article{blatt2009,
  author     = {Blatt, Simon},
  title      = {A singular example for the {W}illmore flow},
  journal    = {Analysis (Munich)},
  fjournal   = {Analysis. International Mathematical Journal of Analysis and
                its Applications},
  volume     = {29},
  year       = {2009},
  number     = {4},
  pages      = {407--430},
  issn       = {0174-4747},
  mrclass    = {53C44 (49Q20)},
  mrnumber   = {2591055},
  mrreviewer = {Yu\ Zheng},
  doi        = {10.1524/anly.2009.1017},
  url        = {https://doi.org/10.1524/anly.2009.1017}
}

@PhdThesis{Jachan,
  author = {Jachan, Felix},
  school = {FU Berlin},
  title  = {Area preserving Willmore flow in asymptotically Schwarzschild manifolds},
  year   = {2014},
}

@article {chillfasangovaschaetzle2009,
    AUTHOR = {Chill, Ralph and Fa\v{s}angov\'{a}, Eva and Sch\"{a}tzle,
              Reiner},
     TITLE = {Willmore blowups are never compact},
   JOURNAL = {Duke Math. J.},
  FJOURNAL = {Duke Mathematical Journal},
    VOLUME = {147},
      YEAR = {2009},
    NUMBER = {2},
     PAGES = {345--376},
      ISSN = {0012-7094,1547-7398},
   MRCLASS = {53C44 (35A15 35K55)},
  MRNUMBER = {2495079},
MRREVIEWER = {Federico\ S\'{a}nchez-Bringas},
       DOI = {10.1215/00127094-2009-014},
       URL = {https://doi.org/10.1215/00127094-2009-014},
}

@article{BKS2024,
url = {https://doi.org/10.1515/acv-2022-0056},
title = {Generalized minimizing movements for the varifold {C}anham--{H}elfrich flow},
author = {Katharina Brazda and Martin Kružík and Ulisse Stefanelli},
journal = {Advances in Calculus of Variations},
doi = {doi:10.1515/acv-2022-0056},
year = {2024},
lastchecked = {2024-06-28}
}

@article {McCoyWheeler2016,
    AUTHOR = {McCoy, James and Wheeler, Glen},
     TITLE = {Finite time singularities for the locally constrained
              {W}illmore flow of surfaces},
   JOURNAL = {Comm. Anal. Geom.},
  FJOURNAL = {Communications in Analysis and Geometry},
    VOLUME = {24},
      YEAR = {2016},
    NUMBER = {4},
     PAGES = {843--886},
      ISSN = {1019-8385,1944-9992},
   MRCLASS = {53C44 (58J35)},
  MRNUMBER = {3570419},
MRREVIEWER = {Andrea\ Mondino},
       DOI = {10.4310/CAG.2016.v24.n4.a7},
       URL = {https://doi.org/10.4310/CAG.2016.v24.n4.a7},
}

@article {Blatt2019,
    AUTHOR = {Blatt, Simon},
     TITLE = {A note on singularities in finite time for the {$L^2$}
              gradient flow of the {H}elfrich functional},
   JOURNAL = {J. Evol. Equ.},
  FJOURNAL = {Journal of Evolution Equations},
    VOLUME = {19},
      YEAR = {2019},
    NUMBER = {2},
     PAGES = {463--477},
      ISSN = {1424-3199,1424-3202},
   MRCLASS = {53C44 (58E20)},
  MRNUMBER = {3950698},
MRREVIEWER = {Glen\ E.\ Wheeler},
       DOI = {10.1007/s00028-019-00483-y},
       URL = {https://doi.org/10.1007/s00028-019-00483-y},
}

@article{choksiveneroni2013,
  author     = {Choksi, Rustum and Veneroni, Marco},
  title      = {Global minimizers for the doubly-constrained {H}elfrich energy: the axisymmetric case},
  journal    = {Calc. Var. Partial Differential Equations},
  fjournal   = {Calculus of Variations and Partial Differential Equations},
  volume     = {48},
  year       = {2013},
  number     = {3-4},
  pages      = {337--366},
  issn       = {0944-2669,1432-0835},
  mrclass    = {49Q10 (49J10 49Q20 53C80 58E30)},
  mrnumber   = {3116014},
  mrreviewer = {Davide\ Vittone},
  doi        = {10.1007/s00526-012-0553-9},
  url        = {https://doi.org/10.1007/s00526-012-0553-9}
}

@article {dallacquamullerschatzlespener2020,
    AUTHOR = {Dall'Acqua, Anna and M\"uller, Marius and Sch\"atzle, Reiner
              and Spener, Adrian},
     TITLE = {The {W}illmore flow of tori of revolution},
   JOURNAL = {Anal. PDE},
  FJOURNAL = {Analysis \& PDE},
    VOLUME = {17},
      YEAR = {2024},
    NUMBER = {9},
     PAGES = {3079--3124},
      ISSN = {2157-5045,1948-206X},
   MRCLASS = {53E40 (49Q20 58E30)},
  MRNUMBER = {4818199},
MRREVIEWER = {Hung\ Thanh\ Tran},
       DOI = {10.2140/apde.2024.17.3079},
       URL = {https://doi.org/10.2140/apde.2024.17.3079},
}

@article{DeulingHelfrich,
	author = {Deuling, HJ and Helfrich, W},
	journal = {Biophysical journal},
	number = {8},
	pages = {861--868},
	publisher = {Elsevier},
	title = {Red blood cell shapes as explained on the basis of curvature elasticity},
	volume = {16},
	year = {1976}
}

@book {Evans98AMS,
    AUTHOR = {Evans, Lawrence C.},
     TITLE = {Partial differential equations},
    SERIES = {Graduate Studies in Mathematics},
    VOLUME = {19},
 PUBLISHER = {American Mathematical Society, Providence, RI},
      YEAR = {1998},
     PAGES = {xviii+662},
      ISBN = {0-8218-0772-2},
   MRCLASS = {35-01},
  MRNUMBER = {1625845},
MRREVIEWER = {Luigi\ Rodino},
       DOI = {10.1090/gsm/019},
       URL = {https://doi.org/10.1090/gsm/019},
}

@article{kuwertschaetzle2001,
  author     = {Kuwert, Ernst and Sch\"{a}tzle, Reiner},
  title      = {The {W}illmore flow with small initial energy},
  journal    = {J. Differential Geom.},
  fjournal   = {Journal of Differential Geometry},
  volume     = {57},
  year       = {2001},
  number     = {3},
  pages      = {409--441},
  issn       = {0022-040X,1945-743X},
  mrclass    = {53C44},
  mrnumber   = {1882663},
  mrreviewer = {Shu-Cheng\ Chang},
  url        = {http://projecteuclid.org/euclid.jdg/1090348128}
}

@article {DeMatteisManno14,
    AUTHOR = {De Matteis, Giovanni and Manno, Gianni},
     TITLE = {Lie algebra symmetry analysis of the {H}elfrich and {W}illmore
              surface shape equations},
   JOURNAL = {Commun. Pure Appl. Anal.},
  FJOURNAL = {Communications on Pure and Applied Analysis},
    VOLUME = {13},
      YEAR = {2014},
    NUMBER = {1},
     PAGES = {453--481},
      ISSN = {1534-0392,1553-5258},
   MRCLASS = {35A30 (35J61 35R01 53A05)},
  MRNUMBER = {3082570},
       DOI = {10.3934/cpaa.2014.13.453},
       URL = {https://doi.org/10.3934/cpaa.2014.13.453},
}

@article {barrettgarckenuernberg2019,
    AUTHOR = {Barrett, John W. and Garcke, Harald and N\"urnberg, Robert},
     TITLE = {Finite element methods for fourth order axisymmetric geometric
              evolution equations},
   JOURNAL = {J. Comput. Phys.},
  FJOURNAL = {Journal of Computational Physics},
    VOLUME = {376},
      YEAR = {2019},
     PAGES = {733--766},
      ISSN = {0021-9991,1090-2716},
   MRCLASS = {65M60 (35Q35)},
  MRNUMBER = {3875543},
MRREVIEWER = {M.\ D.\ Marcozzi},
       DOI = {10.1016/j.jcp.2018.10.006},
       URL = {https://doi.org/10.1016/j.jcp.2018.10.006},
}

@article {mayersimonett2002,
    AUTHOR = {Mayer, Uwe F. and Simonett, Gieri},
     TITLE = {A numerical scheme for axisymmetric solutions of
              curvature-driven free boundary problems, with applications to
              the {W}illmore flow},
   JOURNAL = {Interfaces Free Bound.},
  FJOURNAL = {Interfaces and Free Boundaries. Modelling, Analysis and
              Computation},
    VOLUME = {4},
      YEAR = {2002},
    NUMBER = {1},
     PAGES = {89--109},
      ISSN = {1463-9963,1463-9971},
   MRCLASS = {53C44 (65M99)},
  MRNUMBER = {1877537},
MRREVIEWER = {Emmanuel\ Hebey},
       DOI = {10.4171/IFB/54},
       URL = {https://doi.org/10.4171/IFB/54},
}

@article {barrettgarckenuernberg2008,
    AUTHOR = {Barrett, John W. and Garcke, Harald and N\"urnberg, Robert},
     TITLE = {Parametric approximation of {W}illmore flow and related
              geometric evolution equations},
   JOURNAL = {SIAM J. Sci. Comput.},
  FJOURNAL = {SIAM Journal on Scientific Computing},
    VOLUME = {31},
      YEAR = {2008},
    NUMBER = {1},
     PAGES = {225--253},
      ISSN = {1064-8275,1095-7197},
   MRCLASS = {53C44 (35K55 65M99)},
  MRNUMBER = {2460777},
       DOI = {10.1137/070700231},
       URL = {https://doi.org/10.1137/070700231},
}

@article {barrettgarckenuernberg2016,
    AUTHOR = {Barrett, John W. and Garcke, Harald and N\"urnberg, Robert},
     TITLE = {Computational parametric {W}illmore flow with spontaneous
              curvature and area difference elasticity effects},
   JOURNAL = {SIAM J. Numer. Anal.},
  FJOURNAL = {SIAM Journal on Numerical Analysis},
    VOLUME = {54},
      YEAR = {2016},
    NUMBER = {3},
     PAGES = {1732--1762},
      ISSN = {0036-1429,1095-7170},
   MRCLASS = {65M60 (35K55 53C44 65M12)},
  MRNUMBER = {3508988},
MRREVIEWER = {Pedro\ Gonz\'alez-Casanova},
       DOI = {10.1137/16M1065379},
       URL = {https://doi.org/10.1137/16M1065379},
}

@article{kuwertschaetzle2002,
  author     = {Kuwert, Ernst and Sch\"{a}tzle, Reiner},
  title      = {Gradient flow for the {W}illmore functional},
  journal    = {Comm. Anal. Geom.},
  fjournal   = {Communications in Analysis and Geometry},
  volume     = {10},
  year       = {2002},
  number     = {2},
  pages      = {307--339},
  issn       = {1019-8385,1944-9992},
  mrclass    = {53C44 (35K10)},
  mrnumber   = {1900754},
  mrreviewer = {Shu-Cheng\ Chang},
  doi        = {10.4310/CAG.2002.v10.n2.a4},
  url        = {https://doi.org/10.4310/CAG.2002.v10.n2.a4}
}

@article{kuwertschaetzle2004,
  author     = {Kuwert, Ernst and Sch\"{a}tzle, Reiner},
  title      = {Removability of point singularities of {W}illmore surfaces},
  journal    = {Ann. of Math. (2)},
  fjournal   = {Annals of Mathematics. Second Series},
  volume     = {160},
  year       = {2004},
  number     = {1},
  pages      = {315--357},
  issn       = {0003-486X,1939-8980},
  mrclass    = {58E12 (53C44)},
  mrnumber   = {2119722},
  mrreviewer = {Shu-Cheng\ Chang},
  doi        = {10.4007/annals.2004.160.315},
  url        = {https://doi.org/10.4007/annals.2004.160.315}
}

@article{langersinger1984,
  author     = {Langer, Joel and Singer, David A.},
  title      = {The total squared curvature of closed curves},
  journal    = {J. Differential Geom.},
  fjournal   = {Journal of Differential Geometry},
  volume     = {20},
  year       = {1984},
  number     = {1},
  pages      = {1--22},
  issn       = {0022-040X,1945-743X},
  mrclass    = {58E10 (53C22)},
  mrnumber   = {772124},
  mrreviewer = {Gudlaugur\ Thorbergsson},
  url        = {http://projecteuclid.org/euclid.jdg/1214438990}
}

@article{liyau1982,
  author     = {Li, Peter and Yau, Shing Tung},
  title      = {A new conformal invariant and its applications to the
                {W}illmore conjecture and the first eigenvalue of compact
                surfaces},
  journal    = {Invent. Math.},
  fjournal   = {Inventiones Mathematicae},
  volume     = {69},
  year       = {1982},
  number     = {2},
  pages      = {269--291},
  issn       = {0020-9910,1432-1297},
  mrclass    = {53C40 (53C21 58G25)},
  mrnumber   = {674407},
  mrreviewer = {Yu.\ Burago},
  doi        = {10.1007/BF01399507},
  url        = {https://doi.org/10.1007/BF01399507}
}

@article {MondinoScharrer20,
    AUTHOR = {Mondino, Andrea and Scharrer, Christian},
     TITLE = {Existence and regularity of spheres minimising the
              {C}anham-{H}elfrich energy},
   JOURNAL = {Arch. Ration. Mech. Anal.},
  FJOURNAL = {Archive for Rational Mechanics and Analysis},
    VOLUME = {236},
      YEAR = {2020},
    NUMBER = {3},
     PAGES = {1455--1485},
      ISSN = {0003-9527,1432-0673},
   MRCLASS = {49Q10 (49N60 58E12 92C37)},
  MRNUMBER = {4076069},
MRREVIEWER = {Kohji\ Ohtsuka},
       DOI = {10.1007/s00205-020-01497-4},
       URL = {https://doi.org/10.1007/s00205-020-01497-4},
}

@article {Aleksandrov1962,
    AUTHOR = {Aleksandrov, A. D.},
     TITLE = {Uniqueness theorems for surfaces in the large. {V}},
   JOURNAL = {Amer. Math. Soc. Transl. (2)},
  FJOURNAL = {Amer. Math. Soc. Transl. (2)},
    VOLUME = {21},
      YEAR = {1962},
     PAGES = {412--416},
   MRCLASS = {53.75},
  MRNUMBER = {150710},
}

@article {Lengeler2018,
    AUTHOR = {Lengeler, Daniel},
     TITLE = {Asymptotic stability of local {H}elfrich minimizers},
   JOURNAL = {Interfaces Free Bound.},
  FJOURNAL = {Interfaces and Free Boundaries. Mathematical Analysis,
              Computation and Applications},
    VOLUME = {20},
      YEAR = {2018},
    NUMBER = {4},
     PAGES = {533--550},
      ISSN = {1463-9963,1463-9971},
   MRCLASS = {35Q35 (35B35 35B40 35Q92 76D27)},
  MRNUMBER = {3893418},
MRREVIEWER = {Atanas\ G.\ Stefanov},
       DOI = {10.4171/IFB/411},
       URL = {https://doi.org/10.4171/IFB/411},
}

@article {Rupp2020,
    AUTHOR = {Rupp, Fabian},
     TITLE = {On the {{\L}}ojasiewicz-{S}imon gradient inequality on
              submanifolds},
   JOURNAL = {J. Funct. Anal.},
  FJOURNAL = {Journal of Functional Analysis},
    VOLUME = {279},
      YEAR = {2020},
    NUMBER = {8},
     PAGES = {108708, 33},
      ISSN = {0022-1236,1096-0783},
   MRCLASS = {26D10 (26E15 37B35 46T05)},
  MRNUMBER = {4129083},
MRREVIEWER = {Tom\'{a}\v{s}\ B\'{a}rta},
       DOI = {10.1016/j.jfa.2020.108708},
       URL = {https://doi.org/10.1016/j.jfa.2020.108708},
}

@article{muellerspener2020,
  author     = {M\"{u}ller, Marius and Spener, Adrian},
  title      = {On the convergence of the elastic flow in the hyperbolic
                plane},
  journal    = {Geom. Flows},
  fjournal   = {Geometric Flows},
  volume     = {5},
  year       = {2020},
  number     = {1},
  pages      = {40--77},
  issn       = {2353-3382},
  mrclass    = {53E99 (49K05 65K10)},
  mrnumber   = {4080255},
  mrreviewer = {Rongpei\ Huang},
  doi        = {10.1515/geofl-2020-0002},
  url        = {https://doi.org/10.1515/geofl-2020-0002}
}

@article {ruppscharrer2023,
    AUTHOR = {Rupp, Fabian and Scharrer, Christian},
     TITLE = {Li-{Y}au inequalities for the {H}elfrich functional and
              applications},
   JOURNAL = {Calc. Var. Partial Differential Equations},
  FJOURNAL = {Calculus of Variations and Partial Differential Equations},
    VOLUME = {62},
      YEAR = {2023},
    NUMBER = {2},
     PAGES = {Paper No. 45, 43},
      ISSN = {0944-2669,1432-0835},
   MRCLASS = {53C42 (49Q10 92C10)},
  MRNUMBER = {4525726},
MRREVIEWER = {Luca\ Lussardi},
       DOI = {10.1007/s00526-022-02381-7},
       URL = {https://doi.org/10.1007/s00526-022-02381-7},
}

@article{RandBurton,
	title = {Mechanical Properties of the Red Cell Membrane: I. Membrane Stiffness and Intracellular Pressure},
	journal = {Biophysical Journal},
	volume = {4},
	number = {2},
	pages = {115-135},
	year = {1964},
	issn = {0006-3495},
	doi = {https://doi.org/10.1016/S0006-3495(64)86773-4},
	url = {https://www.sciencedirect.com/science/article/pii/S0006349564867734},
	author = {R.P. Rand and A.C. Burton},
}

@article{rupp2024,
    AUTHOR = {Rupp, Fabian},
     TITLE = {The {W}illmore flow with prescribed isoperimetric ratio},
   JOURNAL = {Comm. Partial Differential Equations},
  FJOURNAL = {Communications in Partial Differential Equations},
    VOLUME = {49},
      YEAR = {2024},
    NUMBER = {1-2},
     PAGES = {148--184},
      ISSN = {0360-5302,1532-4133},
   MRCLASS = {53E40 (35B40 35K41 49Q10)},
  MRNUMBER = {4701428},
       DOI = {10.1080/03605302.2024.2302377},
       URL = {https://doi.org/10.1080/03605302.2024.2302377},
}

@article{rupp2023,
  author     = {Rupp, Fabian},
  title      = {The volume-preserving {W}illmore flow},
  journal    = {Nonlinear Anal.},
  fjournal   = {Nonlinear Analysis. Theory, Methods \& Applications. An
                International Multidisciplinary Journal},
  volume     = {230},
  year       = {2023},
  pages      = {Paper No. 113220, 30},
  issn       = {0362-546X,1873-5215},
  mrclass    = {53E40 (35B40 35K41)},
  mrnumber   = {4541420},
  mrreviewer = {Glen\ E.\ Wheeler},
  doi        = {10.1016/j.na.2023.113220},
  url        = {https://doi.org/10.1016/j.na.2023.113220}
}

@article{canham1970,
  title={The minimum energy of bending as a possible explanation of the biconcave shape of the human red blood cell},
  author={Canham, Peter},
  journal={J. Theor. Biol.},
  fjournal={Journal of theoretical biology},
  volume={26},
  number={1},
  pages={61--81},
  year={1970},
  publisher={Elsevier}
}

@article {bryant1984,
    AUTHOR = {Bryant, Robert},
     TITLE = {A duality theorem for {W}illmore surfaces},
   JOURNAL = {J. Differential Geom.},
  FJOURNAL = {Journal of Differential Geometry},
    VOLUME = {20},
      YEAR = {1984},
    NUMBER = {1},
     PAGES = {23--53},
      ISSN = {0022-040X,1945-743X},
   MRCLASS = {58E12 (53C42)},
  MRNUMBER = {772125},
MRREVIEWER = {Akito\ Futaki},
       URL = {http://projecteuclid.org/euclid.jdg/1214438991},
}

@article{helfrich1973,
  title={Elastic properties of lipid bilayers: theory and possible experiments},
  author={Helfrich, Wolfgang},
  journal={Zeitschrift f{\"u}r Naturforschung C},
  volume={28},
  number={11-12},
  pages={693--703},
  year={1973},
  publisher={Verlag der Zeitschrift f{\"u}r Naturforschung}
}

@article {Scharrer22NLA,
    AUTHOR = {Scharrer, Christian},
     TITLE = {Embedded {D}elaunay tori and their {W}illmore energy},
   JOURNAL = {Nonlinear Anal.},
  FJOURNAL = {Nonlinear Analysis. Theory, Methods \& Applications. An
              International Multidisciplinary Journal},
    VOLUME = {223},
      YEAR = {2022},
     PAGES = {Paper No. 113010, 23},
      ISSN = {0362-546X,1873-5215},
   MRCLASS = {53A10 (33E05)},
  MRNUMBER = {4438233},
MRREVIEWER = {Paul\ F.\ Bracken},
       DOI = {10.1016/j.na.2022.113010},
       URL = {https://doi.org/10.1016/j.na.2022.113010},
}

@article {Schygulla,
    AUTHOR = {Schygulla, Johannes},
     TITLE = {Willmore minimizers with prescribed isoperimetric ratio},
   JOURNAL = {Arch. Ration. Mech. Anal.},
  FJOURNAL = {Archive for Rational Mechanics and Analysis},
    VOLUME = {203},
      YEAR = {2012},
    NUMBER = {3},
     PAGES = {901--941},
      ISSN = {0003-9527,1432-0673},
   MRCLASS = {53C42 (49J10 52A40 74L15 92C37)},
  MRNUMBER = {2928137},
MRREVIEWER = {Luigi\ De Pascale},
       DOI = {10.1007/s00205-011-0465-4},
       URL = {https://doi.org/10.1007/s00205-011-0465-4},
}

@article{simon1993,
  author     = {Simon, Leon},
  title      = {Existence of surfaces minimizing the {W}illmore functional},
  journal    = {Comm. Anal. Geom.},
  fjournal   = {Communications in Analysis and Geometry},
  volume     = {1},
  year       = {1993},
  number     = {2},
  pages      = {281--326},
  issn       = {1019-8385,1944-9992},
  mrclass    = {58E12 (49Q10 53A10)},
  mrnumber   = {1243525},
  mrreviewer = {J.\ E.\ Brothers},
  doi        = {10.4310/CAG.1993.v1.n2.a4},
  url        = {https://doi.org/10.4310/CAG.1993.v1.n2.a4}
}

@article {Willmore65,
    AUTHOR = {Willmore, T. J.},
     TITLE = {Note on embedded surfaces},
   JOURNAL = {An. \c{S}ti. Univ. ``Al. I. Cuza'' Ia\c{s}i Sec\c{t}. I a Mat.
              (N.S.)},
  FJOURNAL = {Analele \c{S}tiin\c{t}ifice ale Universit\u{a}\c{t}ii "Al. I.
              Cuza" din Iasi. Sec\c{t}iunea I a. Matematic\u{a}. Serie
              Nou\u{a}},
    VOLUME = {11B},
      YEAR = {1965},
     PAGES = {493--496},
      ISSN = {0041-9109},
   MRCLASS = {53.01 (53.74)},
  MRNUMBER = {202066},
MRREVIEWER = {H.\ B.\ Griffiths},
}

@book {Willmore93Oxford,
    AUTHOR = {Willmore, T. J.},
     TITLE = {Riemannian geometry},
    SERIES = {Oxford Science Publications},
 PUBLISHER = {The Clarendon Press, Oxford University Press, New York},
      YEAR = {1993},
     PAGES = {xii+318},
      ISBN = {0-19-853253-9},
   MRCLASS = {53-02 (53-01 53C42)},
  MRNUMBER = {1261641},
MRREVIEWER = {Paul\ E.\ Ehrlich},
}

\end{document}